\documentclass{amsart}%
\usepackage{graphicx}
\PassOptionsToPackage{normalem}{ulem}
\usepackage{ulem}

\usepackage{color}

\makeatletter
\usepackage{amsmath}
\usepackage{amstext}
\usepackage{amssymb}

\@ifundefined{showcaptionsetup}{}{%
 \PassOptionsToPackage{caption=false}{subfig}}
\usepackage{subfig}

\newtheorem{theorem}{Theorem}
\newtheorem{lemma}{Lemma}[section]

\newtheorem{proposition}[lemma]{Proposition}

\newtheorem{remark}[lemma]{Remark}

\begin{document}
\title{Delay-induced blow-up in a planar oscillation model}


\author[A. Eremin]{Alexey Eremin} 
\address[A. Eremin]{Department of Information Systems, Saint Petersburg State University} 
\author[E. Ishiwata]{Emiko Ishiwata} 
\address[E. Ishiwata]{Department of Applied Mathematics, Tokyo University of Science} 
\author[T. Ishiwata]{Tetsuya Ishiwata} 
\address[T. Ishiwata]{Department of Mathematical Sciences, Shibaura Institute of Technology}
\email{tisiwata@shibaura-it.ac.jp}
\author[Y. Nakata]{Yukihiko Nakata}
\address[Y. Nakata]{Department of Mathematical Sciences, Aoyama Gakuin University} 


\keywords{blow-up of solutions, periodic solutions, delay differential equations}
\thanks{This work was supported by KAKENHI No. 26400212, No.15K13461, No. 19H05599,  No. 19K21836 and No. 20K03734.}
\subjclass[2010]{34K99, 34K13}

\maketitle

\begin{abstract}
In this paper we study a system of delay differential equations from
the viewpoint of a finite time blow-up of the solution. We prove that
the system admits blow-up solutions, no matter how small the length
of the delay is. In the non-delay system every solution approaches
to a stable unit circle in the plane, thus time delay induces 
blow-up of solutions, 
which we call ``delay-induced blow-up'' phenomenon. Furthermore,
it is shown that the system has a family of infinitely many periodic
solutions, while the non-delay system has only one stable limit cycle.
The system studied in this paper is an example that arbitrary small
delay can be responsible for a drastic change of the dynamics. We
show numerical examples to illustrate our theoretical results.
\end{abstract}

\section{Introduction}

In various disciplines of the science, mathematical modelling offers
a description of phenomena. In some phenomena, the history of the
state, not only the current state, affects the change of the state,
thus it is reasonable to consider the effect of time delay \cite{Erneux:2009}.
Up to now the theory of delay differential equations have been intensively
developed \cite{Hale:1993,Smith:2011}.

In this paper we study a system of delay differential equations from
the viewpoint of the blow-up solutions. Here we use the terminology
``blow-up'' as a finite time blow-up of the solutions, i.e., the solution
diverges (in a suitable topology) in finite time. The blow-up phenomenon
has been widely investigated in partial differential equations, see
\cite{Bandle:1998,Hirota:2006,Straughan:2019,Souplet:1998} and references
therein. There are also extensive studies in Volterra integral equation,
as an alternative formulation of partial differential equation of
parabolic type with a point source term \cite{Brunner:2012,Souplet:1998,Yang:2013}.
Compactification of the phase space is a method to study the blow-up
solutions of ordinary differential equations \cite{Elias:2006,Matsue:2018}.
Numerical analysis has been an unavoidable tool for understanding
the blow-up phenomenon. Numerical method to compute the blow-up solutions
for polynomial systems of ordinary differential equations has been
proposed for ordinary differential equations with application to partial
differential equations, see \cite{Hirota:2006} and references therein.
Recently, numerical validation for the existence of the blow-up solutions
is proposed \cite{Takayasu:2017}.

The blow-up phenomenon in delay differential equations has not been
much studied, as far as we know, except for a few studies \cite{EzzinbiJazar2006,Ma:2011,Yang:2017,Zhou:2016,Hirata:1999,Appleby:2017}.
Perhaps the reason is that
many examples of delay differential equations, which appear in population
biology, control theory, etc, negative feedback condition is usually
imposed which excludes blow-up solutions. One also sees that the following
delay differential equation 
\[
x^{\prime}(t)=x(t)x(t-\tau)
\]
does not have a blow-up solution (at least if the initial function
is continuous) for $\tau>0$ and is reduced to a famous example of
the ordinary differential equation having a blow-up solution when
$\tau=0$. Thus one may speculate that time delay inhibits the blow-up
phenomenon in general (which is certainly not).

In the paper \cite{EzzinbiJazar2006} the authors study the existence
of the blow-up solutions for a class of delay differential equations.
The authors are interested if adding (and multiplying) a delay term
to an ordinary differential equation affects the solution behavior.
Using the comparison principle, the authors obtain conditions that
the delay term does not change the qualitative properties concerning
 the global existence
and blow-up of the solution. In \cite{Zhou:2016} the authors study
the blow-up phenomenon of differential equations with piecewise constant
arguments, in comparison with the corresponding ordinary differential
equations. The blow-up phenomenon is studied in Volterra integro-differential
equations. See \cite{Hirata:1999,Appleby:2017} and also
\cite{Ma:2011,Yang:2017} as applications
to a parabolic type partial differential equations to study
of Volterra integro-differential equations.

It seems that the blow-up phenomenon that stems from the time
delay has not been reported, to the best of our knowledge. 
Our motivations are to demonstrate whether the time delay itself
 induces blow-up of the solution, and to understand the mechanism of blow-up of solutions. 
In this paper, we propose an example model that time delay drastically changes the
solution behavior and induces the blow-up solution together with a
family of infinitely many periodic solutions, where most of solutions
are shown to be unstable. In our equation, a blow-up solution exists
\textit{no matter how small the length of the delay is}, while non-delay
equation does not have a blow-up solution. We thus call these phenomena
``delay-induced blow-up''.

This paper is organized as follows. In Section 2, we introduce a planar
system of delay differential equations which we study in this paper.
In the absence of time delay, the system becomes a planar system of
ordinary differential equations. It can be seen that every solution
except the trivial solution approaches to the limit cycle, thus no
solution blows up. Concerning the system of delay differential equations, we
present our two main theorems that show blow-up of solutions is
possible due to time delay and that the system admits infinitely many
periodic solutions.
In this section we use a special initial function in order to show
 blow-up of solutions for any $\tau>0$. We demonstrate numerical examples
 of blow-up solutions and global solutions. See Remark \ref{rmk:numerics}, 
Figure \ref{fig:numex} and also Section \ref{Numerics}.
In Section 3, by a careful estimation of the solution we show that 
there is a blow-up solution for the system of delay differential 
equations and provide a proof of Theorem \ref{thm:BU}. 
In Section 4, we study existence of a periodic solution with constant radius
and constant angular velocity in the plane. It is shown that there exist infinitely
many periodic solutions which appear due to time delay. 
In Section 5, we study a characteristic equation which characterizes stability
of the periodic solutions. 
It is shown that most of periodic solutions are not stable except for 
the only one periodic solution that is a continuation of the periodic solution 
of the non-delay model. 
In Section 6, we provide numerical examples which illustrate our theoretical results.
In Section 7, we discuss our results.

\section{A planar system of delay differential equations and main results}

Let $\tau\geq0$ be a parameter for time delay. In this paper we consider
the following planar system of delay differential equations \begin{subequations}\label{eq:DDE}
\begin{align}
x^{\prime}(t) & =x(t)-y(t)-x(t-\tau)\left(x^{2}(t)+y^{2}(t)\right),\label{eq:DDE1}\\
y^{\prime}(t) & =x(t)+y(t)-y(t-\tau)\left(x^{2}(t)+y^{2}(t)\right).\label{eq:DDE2}
\end{align}
\end{subequations}For the special case $\tau=0$, the system (\ref{eq:DDE})
is reduced to the following system of ordinary differential equations 
with $a=1$:
\begin{subequations}\label{eq:ODE}
\begin{align}
x^{\prime}(t) & =a x(t)-y(t)-x(t)\left(x^{2}(t)+y^{2}(t)\right),\label{eq:ODE1}\\
y^{\prime}(t) & =x(t)+a y(t)-y(t)\left(x^{2}(t)+y^{2}(t)\right).\label{eq:ODE2}
\end{align}
\end{subequations}
Here $a\in \mathbb{R}$ is a parameter.
This system is a famous model for Hopf bifurcation at $a=0$. 
(See \cite{Kuznetsov:1998}, for instance.)
From the elementary calculation, one can see that
every solution of (\ref{eq:ODE}) except the trivial solution $\left(x\left(t\right),y\left(t\right)\right)\equiv\left(0,0\right)$
tends to a periodic solution of minimal period $2\pi$ and satisfies
\[
\lim_{t\to\infty}\sqrt{x^{2}(t)+y^{2}(t)}=1.
\]
i.e., the limit cycle of (\ref{eq:ODE}) is the unit circle.

For the system (\ref{eq:DDE}) we prove that the system admits blow-up solutions
 due to the presence of time delay.  We prove the following
theorem in Section 3.

\begin{theorem}
\label{thm:BU}
For any $\tau>0$, 
there exist blow-up solutions 
for the planar system of delay differential equations (\ref{eq:DDE}).
\end{theorem}

Then in Sections 4 and 5, we further investigate the system (\ref{eq:DDE}) and
show the existence of infinitely many periodic solutions. We also
study the stability of the periodic solutions. The following theorem
is proved in Sections 4 and 5.

\begin{theorem}\label{thm:periodic}
For any $\tau>0$, there exist infinitely many unstable periodic solutions
of (\ref{eq:DDE}) with constant radius and angular velocity.
Moreover, there exists a positive $\tau^*$ such that
the system (\ref{eq:DDE}) admits only one asymptotically stable periodic solution
 with constant radius and angular velocity 
 for $0<\tau<\tau^*$.
\end{theorem}

\section{Delay-induced blow-up}

We consider a polar coordinate system. Let
\[
(x(t),y(t))=r(t)\left(\cos\theta(t),\sin\theta(t)\right).
\]
By the change of the variables, from (\ref{eq:DDE}), we obtain the
following polar coordinate system\begin{subequations}\label{eq:DDE_PC}
\begin{align}
r^{\prime}(t) & =r(t)\left(1-r(t)r(t-\tau)\cos\left(\theta(t)-\theta(t-\tau)\right)\right),\label{eq:DDE1_PC}\\
\theta^{\prime}(t) & =1+r(t)r(t-\tau)\sin\left(\theta(t)-\theta(t-\tau)\right).\label{eq:DDE2_PC}
\end{align}
\end{subequations}Now we prove that there exists a solution such
that $r$ blows up in a finite time with $\theta\to\frac{\pi}{4}$
as $r\to\infty$. Therefore, we can conclude that $x$ and $y$ blow up
in a finite time.

For (\ref{eq:DDE_PC}), we consider the following initial condition\begin{subequations}\label{eq:IC}
\begin{align}
r(t) & =R,\ t\in\left[-\tau,0\right],\label{eq:PCIC1}\\
\theta(t) & =\begin{cases}
-\frac{3}{4}\pi, & -\tau\leq t<-\frac{1}{2}\tau,\\
s(t), & -\frac{1}{2}\tau\leq t\leq0,
\end{cases}\label{eq:PCIC2}
\end{align}
\end{subequations}where $R>0$ and $s(t)$ is a continuous function
such that 
\[
s\left(-\frac{1}{2}\tau\right)=-\frac{3}{4}\pi,\ s\left(0\right)=-\frac{1}{2}\pi.
\]

\begin{figure}
\begin{center}
\begin{tabular}{ccc}
\includegraphics[scale=0.25]{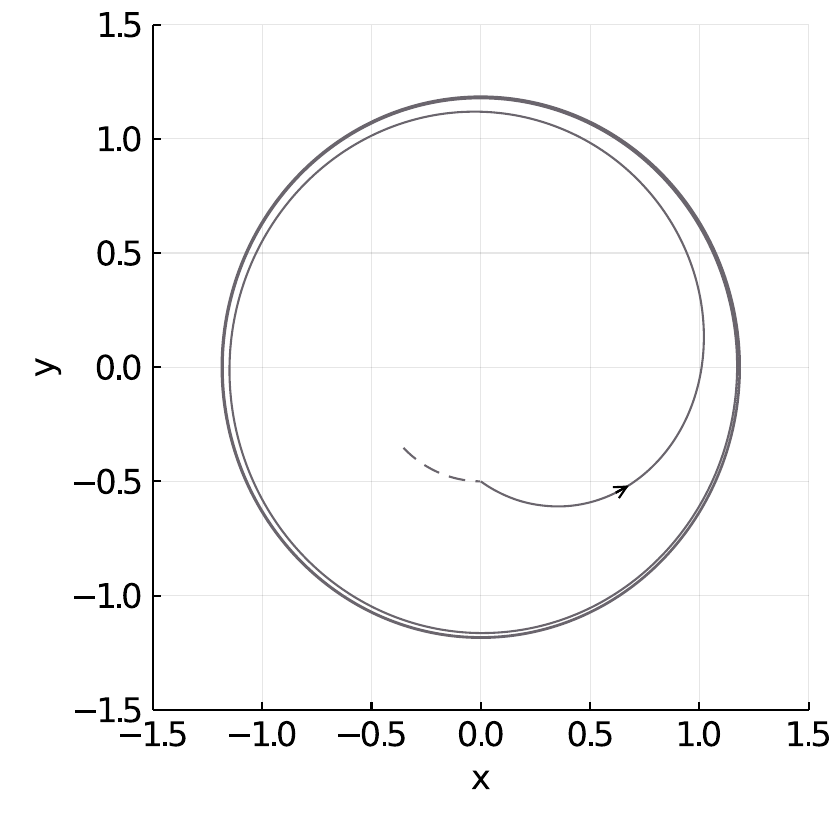} &
\includegraphics[scale=0.25]{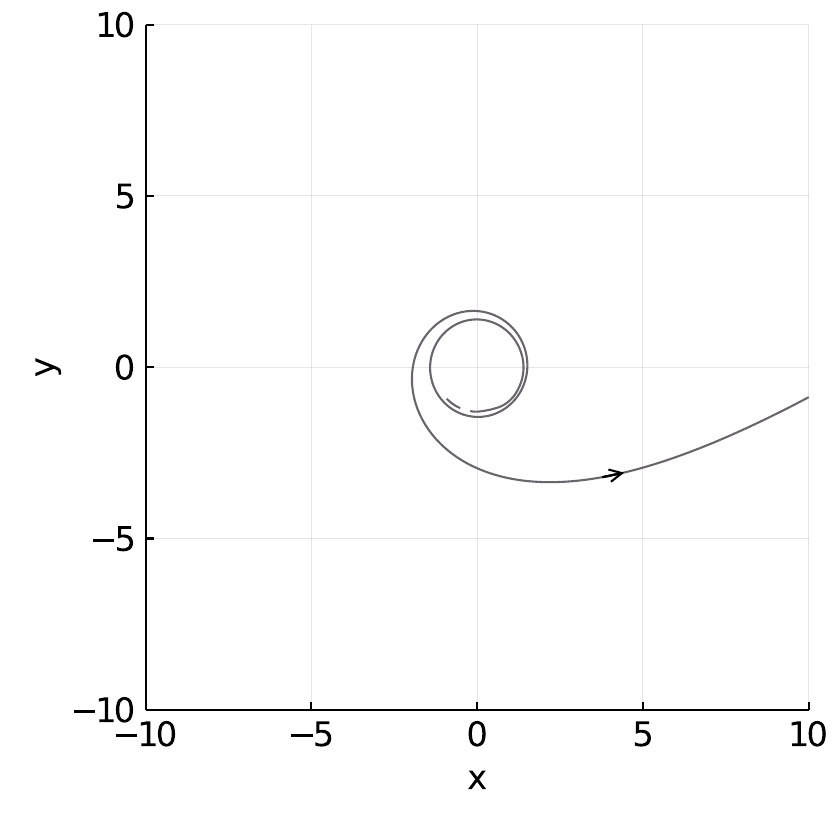} &
\includegraphics[scale=0.25]{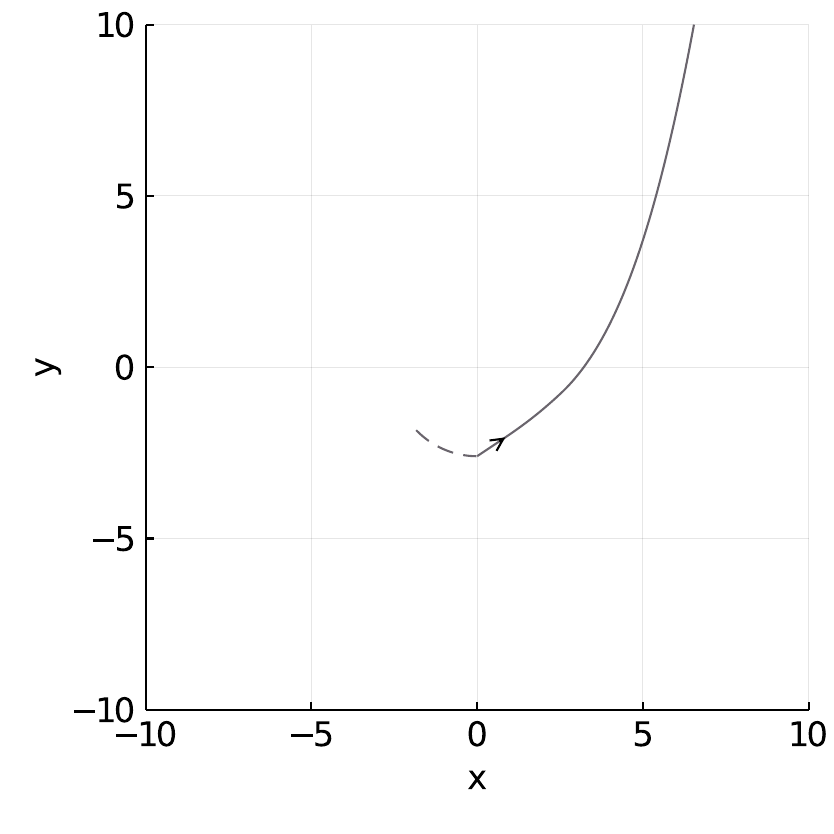} \\
\hspace*{4mm} (i) & \hspace*{4mm} (ii) & \hspace*{4mm}(iii) 
\end{tabular}
\caption{\label{fig:numex}
Illustration of trajectories of the solution in $(x,y)$-plane for $\tau=0.392$. 
The initial condition (\ref{eq:IC}), which are plotted as dashed curves 
in the above figures, 
 is used (in the polar coordinate) with 
$s(t)=\frac{\pi}{2\tau}t-\frac{1}{2}\pi$ for the numerical experiments. 
The solution behavior changes by different $R=0.5, 1.3$ and $2.6$. 
(i) The solution converges to a periodic solution with $R=0.5$, 
(ii) The solution blows up in a finite time  $t>\tau$ with $R=1.3$ 
(The blow-up time is nearly $t=4.359\cdots$.)
and 
(iii) The solution blows up in $t<\tau$ with $R=2.6$.
}
\end{center}
\end{figure}

\begin{remark}\label{rmk:numerics}
In this paper we consider a special initial condition (\ref{eq:IC}) in order to
 show an existence of blow-up solutions for any positive $\tau$.
In this section we will prove the blow-up of solutions in $t<\tau$
for sufficiently large $R>0$.
For other initial data, there are no mathematical results on 
global existence of solutions except periodic solutions studied in Section
4 and blow-up of the solutions. 
In Figure \ref{fig:numex} (i) and (ii), numerical examples for 
$R=0.5$ and $1.3$ are presented.
These numerical simulations suggest that there are not only time-global solutions
 which converge to a periodic orbit but also blow-up solutions which blow up 
in a finite time $t>\tau$.
In Section \ref{Numerics}, numerical figures for other initial data
are presented.
\end{remark}

Since  $-\tau\leq t-\tau<-\frac{1}{2}\tau$ for $t\in\left[0,\frac{1}{2}\tau\right]$,
from the initial condition (\ref{eq:IC}), we have 
\[
r(t-\tau)=R,\ \theta(t-\tau)=-\frac{3}{4}\pi,\ t\in\left[0,\frac{1}{2}\tau\right].
\]
Then, for $t\in\left[0,\frac{1}{2}\tau\right]$, the solution of the
system of delay differential equations (\ref{eq:DDE_PC}) with the
initial condition (\ref{eq:IC}) is given by the following system
of ordinary differential equations\begin{subequations}\label{eq:ODEinit}
\begin{align}
r^{\prime}(t) & =r(t)\left(1-Rr(t)\cos\left(\theta(t)+\frac{3}{4}\pi\right)\right),\label{eq:ODEinit1}\\
\theta^{\prime}(t) & =1+Rr(t)\sin\left(\theta(t)+\frac{3}{4}\pi\right)\label{eq:ODEinit2}
\end{align}
\end{subequations}with the initial condition
\begin{equation}
\left(r(0),\theta(0)\right)=\left(R,-\frac{1}{2}\pi\right).\label{eq:ODEinitIC}
\end{equation}

We are going to prove that the solution of the system of ordinary differential equations
(\ref{eq:ODEinit}) with the initial condition (\ref{eq:ODEinitIC})
blows up in $t\leq\frac{1}{2}\tau$. 

\begin{figure}
\includegraphics[scale=1.0]{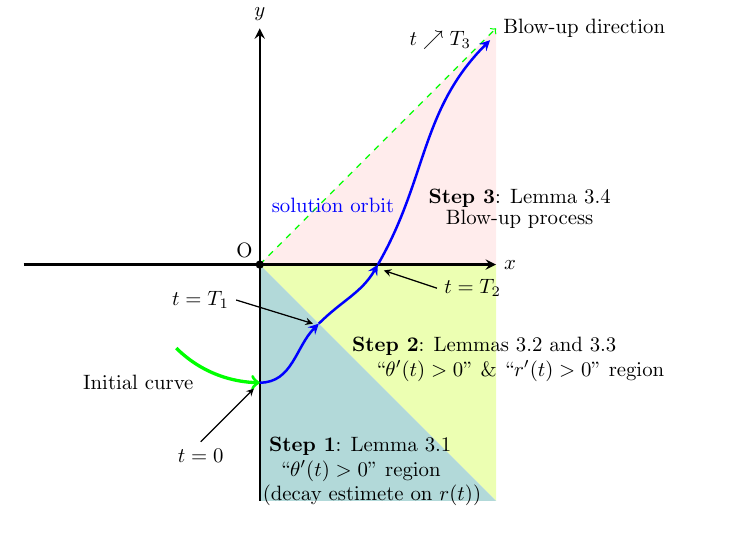} \\
\caption{
Three steps for the blow-up of solutions.
}
\label{fig:steps-for-blowup}
\end{figure}

Figure \ref{fig:steps-for-blowup} shows each step for blow-up of solutions.
The process for blow-up of solutions is divided into the following 3 steps:

\begin{description}
\item[Step 1] The angle of the solutions is monotonically increasing 
from $-\pi/2$ to $-\pi/4$ in $t \in [0, T_1]$ for some $T_1>0$. 
In this region, the radius of the solutions may decrease, and thus
 we establish a decay estimate of the solutions in Lemma \ref{lem:lem1}.

\item[Step 2] Both radius and angle of the solutions monotonically increase
 in $t\in [T_1, T_2]$ for some $T_2>T_1$. The angle varies from $-\pi/4$ to $0$
 and the radius grows up beyond a threshold for an emergence of a nullcline of the 
angle. (See Lemmas \ref{lem:lem2} and \ref{lem:lem34}.)

\item[Step 3] The final stage for blow-up of solutions. The angle monotonically reaches
 to the nullcline of the angle which appears in Step 2.
Then we show in Lemma \ref{lem:lem3} that for large $R$ there exists a blow-up time $T_{3}$
($T_{2}<T_{3})$ with $T_{3}<\frac{1}{2}\tau$ such that 
\begin{equation}
\lim_{t\uparrow T_{3}}r(t)=\infty.
\label{eq:rblowup}
\end{equation}

\end{description}

For the exposition, we define 
\[
\phi(t):=\theta(t)+\frac{3}{4}\pi,\ t\geq0.
\]
Then the system (\ref{eq:ODEinit}) with the initial condition (\ref{eq:ODEinitIC})
becomes \begin{subequations}\label{eq:ODEinitphi}
\begin{align}
r^{\prime}(t) & =r(t)\left(1-Rr(t)\cos\phi(t)\right),\label{eq:ODEinitphi1}\\
\phi^{\prime}(t) & =1+Rr(t)\sin\phi(t)\label{eq:ODEinitphi2}
\end{align}
\end{subequations}with the initial condition
\begin{equation}
\left(r(0),\phi(0)\right)=\left(R,\frac{1}{4}\pi\right).\label{eq:ODEinitICphi}
\end{equation}

In the proof below we need to estimate the solution $r$. For the
estimation we use the following equation 
\begin{equation}
\frac{dr}{d\phi}=\frac{r\left(1-Rr\cos\phi\right)}{1+Rr\sin\phi},\label{eq:drdphi}
\end{equation}
which is obtained from (\ref{eq:ODEinitphi}).

\subsection*{Step 1: The decay estimate of the radius}

\begin{lemma}
\label{lem:lem1}There exists $T_{1}>0$ such that $\phi$ monotonically
increases from $\pi/4$ to $\pi/2$ for $t\in\left[0,T_{1}\right]$.
One also has 
\begin{equation}
R\exp\left(-\left(\phi(t)-\frac{\pi}{4}\right)\right)\leq r(t)\leq R\exp\left(\phi(t)-\frac{\pi}{4}\right),\ t\in\left[0,T_{1}\right].\label{eq:rel_Rtheta}
\end{equation}
Moreover, there exists $\overline{R}_{1}>0$ such that 
\[
R>\overline{R}_{1}\implies T_{1}<\frac{\tau}{4}.
\]
\end{lemma}

\begin{proof}
The solution of the system of ordinary differential equation (\ref{eq:ODEinitphi})
with the initial condition (\ref{eq:ODEinitICphi}) exists for sufficiently
small $t$. 
We show that 
there exists $T_1>0$ such that the solution exists for $t\in\left[0,T_{1}\right]$
satisfying that $\phi$ monotonically increases from $\pi/4$ to $\pi/2$ 
for $0\leq t\leq T_{1}$.
 We first obtain an a priori estimate for $r$
for $\phi\in\left[\frac{1}{4}\pi,\frac{1}{2}\pi\right]$. From (\ref{eq:drdphi})
it follows that 
\begin{equation}
\frac{dr}{d\phi}\leq r,\ \phi\in\left[\frac{1}{4}\pi,\frac{1}{2}\pi\right].\label{eq:drdthtea1}
\end{equation}
One also obtains the following estimation
\begin{equation}
\frac{dr}{d\phi}\geq\frac{-Rr^{2}\cos\phi}{Rr\sin\phi}\geq\frac{-Rr^{2}\frac{\sqrt{2}}{2}}{Rr\frac{\sqrt{2}}{2}}=-r,\ \phi\in\left[\frac{1}{4}\pi,\frac{1}{2}\pi\right].\label{eq:drdtheta2}
\end{equation}
Integrating the inequalities (\ref{eq:drdthtea1}) and (\ref{eq:drdtheta2}),
we obtain the following estimation
\begin{equation}
R\exp\left(-\left(\phi-\frac{\pi}{4}\right)\right)\leq r\leq R\exp\left(\phi-\frac{\pi}{4}\right),\ \phi\in\left[\frac{1}{4}\pi,\frac{1}{2}\pi\right].\label{eq:apriori_r1}
\end{equation}
Using the a priori bound (\ref{eq:apriori_r1}), from the equation
(\ref{eq:ODEinitphi2}), we see that $\phi$ is an increasing function
and that $\phi^{\prime}(t)\geq1$, provided $\phi\in\left[\frac{1}{4}\pi,\frac{1}{2}\pi\right].$
Therefore, there exists $T_{1}>0$ such that $\phi$ monotonically
increases from $\pi/4$ to $\pi/2$ for $t\in\left[0,T_{1}\right]$. The
inequality (\ref{eq:rel_Rtheta}) holds from the estimation (\ref{eq:apriori_r1}).
Since, from the inequality (\ref{eq:rel_Rtheta}), we have 
$R\exp(-\pi/4)\leq r$
for $t\in\left[0,T_{1}\right]$, the following estimation 
\begin{equation}
\phi^{\prime(t)\geq1+R^{2}\exp\left(-\frac{\pi}{4}\right)}\sin\phi\geq
1+R^{2}\exp\left(-\frac{\pi}{4}\right)\frac{\sqrt{2}}{2},\ t\in\left[0,T_{1}\right]\label{eq:theta_est1}
\end{equation}
implies that there exists $\overline{R}_1>0$ such that if $R>\overline{R}_1$
 then $T_{1}< \frac{\tau}{4}$.
\qed
\end{proof}
From Lemma \ref{lem:lem1}, we have 
\begin{align*}
 & R\exp\left(-\frac{\pi}{4}\right)\leq r(T_{1})\leq R\exp\left(\frac{\pi}{4}\right),\\
 & \phi(T_{1})=\frac{1}{2}\pi.
\end{align*}

\subsection*{Step 2: Emergence of the nullcline of the angle}

Next we have the following estimation.
\begin{lemma}
\label{lem:lem2}There exists $T_{2}\ (>T_{1})$ such that $\phi$
monotonically increases from $\pi/2$ to $3\pi/4$ 
for $t\in\left[T_{1},T_{2}\right]$. 
One also has that $r$ monotonically increases for $t\in\left[T_{1},T_{2}\right]$
and 
\begin{equation}
R\exp\left(-\frac{\pi}{4}\right)\leq r(t)<\infty,\ t\in\left[T_{1},T_{2}\right].\label{eq:rel_Rtheta2}
\end{equation}
Moreover, there exists $\overline{R}_{2}>0$ such that 
\[
R>\overline{R}_{2}\implies T_{2}-T_{1}<\frac{\tau}{8}.
\]
\end{lemma}

\begin{proof}
First, we derive an a priori estimate for $r$, provided $\phi\in\left[\frac{1}{2}\pi,\frac{3}{4}\pi\right]$.
Let $\phi\in\left[\frac{1}{2}\pi,\frac{3}{4}\pi\right]$. We have
$r^{\prime}(t)\geq r$ from (\ref{eq:ODEinitphi1}), thus $r$ monotonically
increases. Then one has
\begin{equation}
r\geq r(T_{1})\geq R\exp\left(-\frac{\pi}{4}\right).\label{eq:apriorir2}
\end{equation}
Since, from (\ref{eq:drdphi}), it holds that 
\begin{equation}
\frac{dr}{d\phi}\leq\frac{r\left(1+Rr\frac{\sqrt{2}}{2}\right)}{Rr\frac{\sqrt{2}}{2}}=\frac{\sqrt{2}+Rr}{R},\ \phi\in\left[\frac{1}{2}\pi,\frac{3}{4}\pi\right],\label{eq:drdthtea2}
\end{equation}
integrating the equation, we obtain 
\begin{equation}
\log\left(\frac{\sqrt{2}+Rr}{\sqrt{2}+Rr(T_{1})}\right)\leq\phi-\frac{\pi}{2},\label{eq:apriorir3}
\end{equation}
which implies $r<\infty$.

From the equation (\ref{eq:ODEinitphi2})
and the a priori bounds (\ref{eq:apriorir2}),
we have that 
$\phi^{\prime}(t)\geq1$ holds for $\phi\in\left[\frac{1}{2}\pi,\frac{3}{4}\pi\right].$ 
Hence, 
there exists $T_{2}>T_{1}$ such that $\phi$ monotonically increases
for $t\in\left[T_{1},T_{2}\right]$ with $\phi(T_{1})=\frac{1}{2}\pi,\ \phi(t)\in\left(\frac{1}{2}\pi,\frac{3}{4}\pi\right)$
for $t\in\left(T_{1},T_{2}\right)$ and $\phi(T_{2})=\frac{3}{4}\pi$.
Therefore, 
the lower estimate (\ref{eq:apriorir2}) is valid for $t\in [T_1, T_2]$
and by virtue of (\ref{eq:apriorir3}) the boundedness of $r(t)$ also holds
for $t\in [T_1, T_2]$.
That is, 
the inequality (\ref{eq:rel_Rtheta2}) holds for $t\in [T_1, T_2].$
Finally, from the
inequality (\ref{eq:rel_Rtheta2}), the following estimation holds:
\begin{equation}
\phi^{\prime}(t)\geq1+R^{2}\exp\left(-\frac{\pi}{4}\right)\sin\phi\geq
1+R^{2}\exp\left(-\frac{\pi}{4}\right)\frac{\sqrt{2}}{2},\ t\in\left[T_{1},T_{2}\right],\label{eq:theta_est1-1}
\end{equation}
which implies that there exists $\overline{R}_2>0$ such that if $R>\overline{R}_2$
 then $T_{2}-T_{1}<\frac{\tau}{8}$.
\qed
\end{proof}

We are ready to show that the solution of (\ref{eq:ODEinitphi}) blows up.
Note that $r'(T_2)\geq r(T_2)>0$ and $\phi'(T_2)>1$ since
$\cos\phi(T_{2})  =\cos\frac{3}{4}\pi=-\frac{\sqrt{2}}{2}<0$ and 
$\sin\phi(T_{2})  =\sin\frac{3}{4}\pi=\frac{\sqrt{2}}{2}>0.$
Thus, there is sufficiently small $\delta>0$ such that for
$t\in\left(T_{2},T_{2}+\delta\right)$ the solution exists and that
$\phi$ and $r$ increase and thus 
$r(t)\geq r(T_{2})\geq R\exp\left(-\frac{\pi}{4}\right)$
for $t\in\left(T_{2},T_{2}+\delta\right)$.

\begin{figure}
\begin{center}
\includegraphics[scale=0.4]{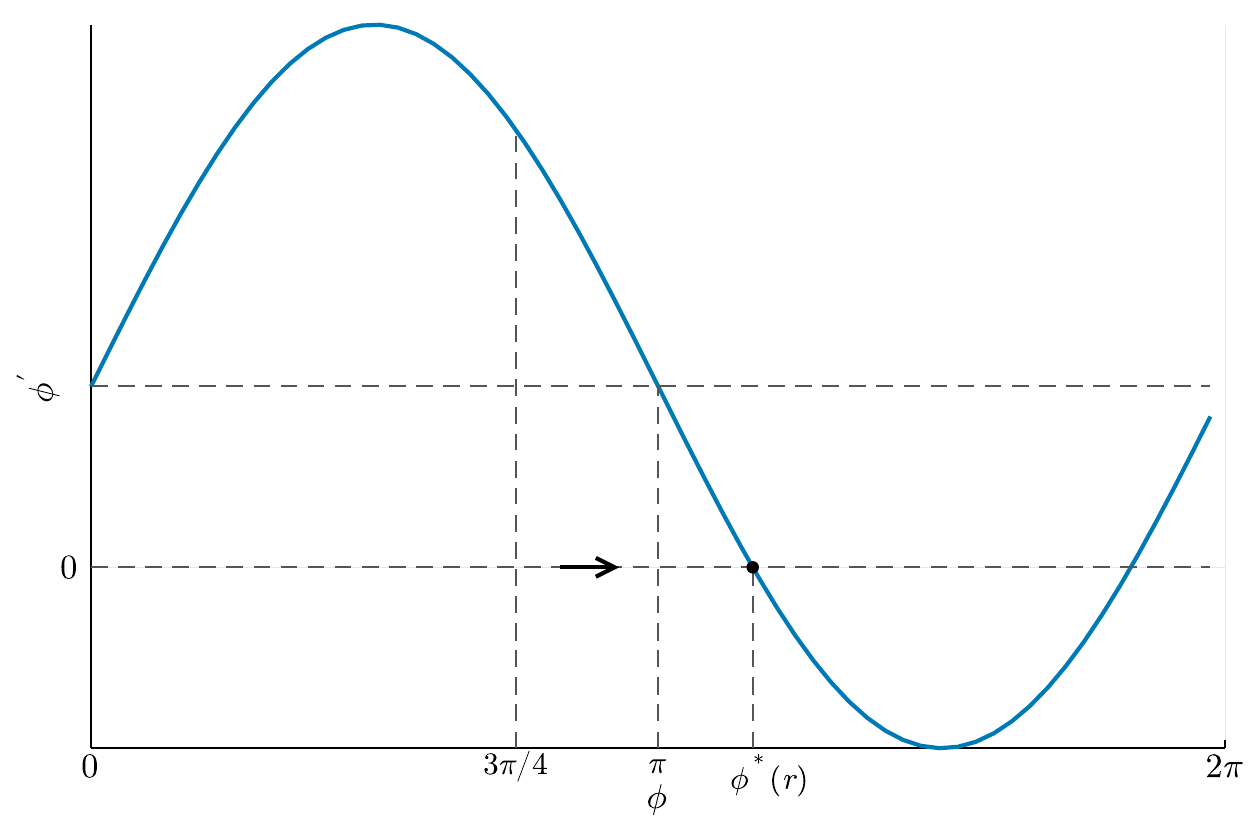}\caption{\label{fig:singraph}
The graph of the right-hand side of the
equation (\ref{eq:ODEinitphi2}) for large $Rr$. Large $Rr$ makes
a nullcline $\phi^{*}(r)$ of the equation (\ref{eq:ODEinitphi2}).
$\phi$ approaches to $\phi^{*}(r)$.}
\end{center}
\end{figure}

In Figure \ref{fig:singraph}, we plot the graph of $\left(\phi,\phi^{\prime}\right)$.
We show that, for sufficiently large $R$, a nullcline for $\phi$
exists in $\left(\pi,\frac{3}{2}\pi\right)$, where the right-hand
side of the $\phi$-equation (\ref{eq:ODEinitphi2}) becomes $0$.
We see in Lemma \ref{lem:lem3} that
$\phi$ has a upper bound for suitable large $R$
and then $r(t)$ blows up in a finite time. 
We now
let $R$ be a sufficiently large number such that $R^{2}\exp(-\frac{\pi}{4})>\sqrt{2}$.
Then, for $r>R\exp(-\frac{\pi}{4})$, there is a $\phi$-nullcline in
$\left(\pi,\frac{5}{4}\pi\right)$ that is given as 
\[
\phi^{*}(r):=\pi-\arcsin\left(-\frac{1}{Rr}\right)\in\left(\pi,\frac{5}{4}\pi\right),\ r>R\exp\left(-\frac{\pi}{4}\right).
\]

It is easy to obtain the following elementary lemma.
\begin{lemma}
\label{lem:lem34}Let $R$ be a sufficiently large number such that
$R^{2}\exp(-\frac{\pi}{4})>\sqrt{2}$. Then
 $\phi^{*}(r)\in\left(\pi,\frac{5}{4}\pi\right)\ 
\left(r>R\exp\left(-\frac{\pi}{4}\right)\right)$.
One has that
\[
1+Rr\sin\phi\begin{cases}
>0, & \phi\in\left(\frac{3}{4}\pi,\phi^{*}(r)\right)\\
=0, & \phi=\phi^{*}(r)
\end{cases}
\]
and that 
\begin{equation}
\lim_{r\to\infty}\phi^{*}(r)=\pi.\label{eq:asymp_phi}
\end{equation}
\end{lemma}

Hence, as long as the solution exists, for $t>T_{2}$, $r(t)\geq r(T_{2})\geq R\exp\left(-\frac{\pi}{4}\right)$
increases and $\phi$ increases and tends to $\phi^{*}(r(t))$. 
Observe that $\phi$ stays in the interval $\left(\frac{3}{4}\pi,\frac{5}{4}\pi\right)$, 
thus
the sign of $\cos\phi$ in the right hand side of (\ref{eq:ODEinitphi1})
is fixed and is positive.

From Lemma \ref{lem:lem2}, we have 
\begin{align*}
 & R\exp\left(-\frac{\pi}{4}\right)\leq r(T_{2})<\infty,\\
 & \phi(T_{2})=\frac{3}{4}\pi.
\end{align*}

\subsection*{Step 3: Blow-up of solutions}

Finally we show the blow-up of solutions.

\begin{lemma}
\label{lem:lem3}There exists $T_{3}\ (>T_{2})$ such that $\phi$
and $r$ monotonically increase for $t\in\left[T_{2},T_{3}\right)$
and 
\begin{equation}
\lim_{t\uparrow T_{3}}r(t)=\infty,\ \lim_{t\uparrow T_{3}}\phi(t)=\pi.\label{eq:asymrphi}
\end{equation}
Furthermore, there exists $\overline{R}_{3}>0$ such that 
\[
R>\overline{R}_{3}\implies T_{3}-T_{2}<\frac{\tau}{8}.
\]
\end{lemma}

\begin{figure}
\begin{center}
\includegraphics[scale=0.4]{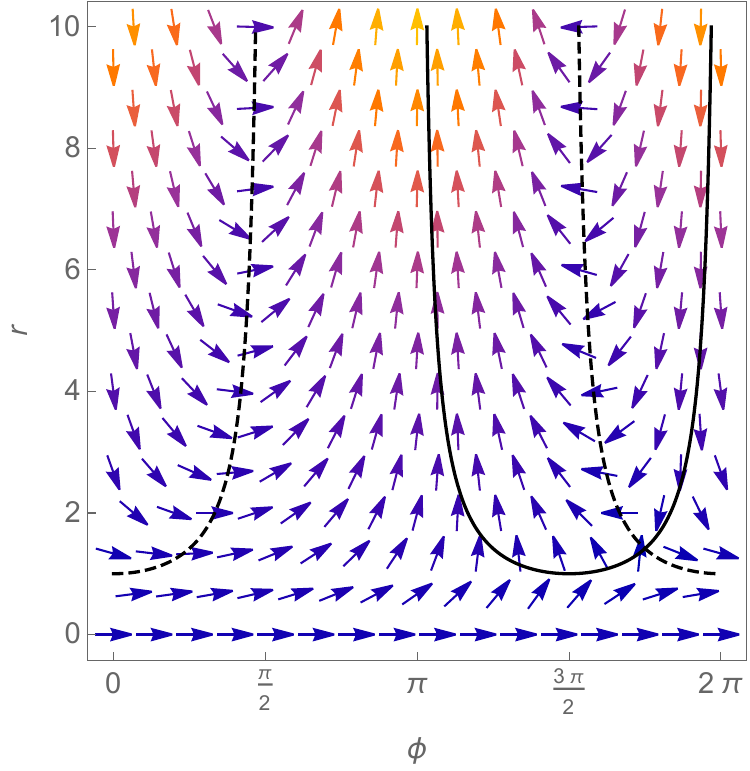}\\
\caption{\label{fig:vectfield}
The vector field $(r'(t), \phi'(t))$ for the system (\ref{eq:ODEinitphi}).
Solid lines (resp. dashed lines) describe $\phi$-nullcline (resp. 
$r$-nullcline).
}
\end{center}
\end{figure}

\begin{proof}
We consider the system (\ref{eq:ODEinitphi}) for $t>T_{2}$. 
Let $R$ be a sufficiently large number such that
 $R>\max\left\{ \overline{R}_{1},\overline{R}_{2}\right\} $
and that $R^{2}\exp(-\frac{\pi}{4})> \sqrt{2}$. 
Since $R>\max\left\{ \overline{R}_{1},\overline{R}_{2}\right\} $,
from Lemmas \ref{lem:lem1} and \ref{lem:lem2}, we have $\phi(T_{2})=\frac{3}{4}\pi$
and $T_{2}<\frac{3}{8}\tau$. We derive an a priori estimate for $\phi$,
provided $r<\infty$. 
Since we have $R^{2}\exp(-\frac{\pi}{4})>\sqrt{2}$,
one sees that the equation (\ref{eq:ODEinitphi2}) has an equilibrium,
$\phi^{*}(r)$ 
$\in (\pi, \frac{5}{4}\pi)$ from Lemma \ref{lem:lem34}. 
Suppose that there exists $t^*>T_2$ such that $\phi'(t^*)=0$ while $r(t^*)< \infty$.
Note that $\phi(t^*)\in (\pi, \frac{5}{4}\pi)$ and by (\ref{eq:drdphi}) we have
$\frac{dr}{d\phi}\vert_{t=t^*}=\infty$, that is, the solution orbit crosses the curve
 of $\phi$-nullcline $\{(\phi, r)| 1+ R r \sin \phi=0, \pi<\phi<\frac{5}{4}\pi \}$
vertically. Then, the solution enters the region where $\phi'(t)<0$, in what follows
$\phi(t)< \frac{5}{4}\pi$. Note that the intersection of the solution orbit and $\phi$-nullcline
 may occur at most once because of the shape of the nullcline and the fact of vertical crossing.
See Figure \ref{fig:vectfield}.

One sees that 
$\phi(t)\in\left(\frac{3}{4}\pi,\frac{5}{4}\pi\right)$.
Since 
\[
-1\leq\cos\phi\leq-\frac{\sqrt{2}}{2},\ \phi\in\left(\frac{3}{4}\pi,\frac{5}{4}\pi\right),
\]
from (\ref{eq:ODEinitphi1}) we now have the estimation,
\begin{align}
r^{\prime}(t) & \geq r\left(1+R\frac{\sqrt{2}}{2}r\right).\label{eq:blowr}
\end{align}
Therefore, one sees that there exists $T_{3}$ such that $\lim_{t\uparrow T_{3}}r(t)=\infty$.
Thus, from (\ref{eq:asymp_phi}) in Lemma \ref{lem:lem34}, we have
$\lim_{t\uparrow T_{3}}\phi(t)=\frac{\pi}{2}$. From the equation (\ref{eq:blowr}),
for any $\varepsilon$ there exists $\overline{R}$ such that $R>\overline{R}$
implies that $T_{3}-T_{2}<\varepsilon$. Therefore, there exists $\overline{R}_{3}$,
such that $T_{3}-T_{2}< \frac{\tau}{8}$. 
\qed
\end{proof}
\begin{flushleft}
\textbf{Proof of Theorem \ref{thm:BU}} Let $R>\overline{R}_{3}$.
From Lemmas \ref{lem:lem1}, \ref{lem:lem2} and \ref{lem:lem3},
one has $T_{3}<\frac{\tau}{2}$. (\ref{eq:asymrphi}) in Lemma \ref{lem:lem3}
implies that
\[
\lim_{t\uparrow T_{3}}\left(x(t),y(t)\right)=\left(\infty,\infty\right).
\]
Thus we obtain the conclusion. \hfill{}$\square$ 
\par\end{flushleft}

\section{Existence of periodic solutions}

We consider a periodic solution with a constant radius,
i.e. $r(t)\equiv \rho>0$, for the system (\ref{eq:DDE_PC}). 
Then, 
$\theta'(t)$ is also a constant 
from 
equations (\ref{eq:DDE_PC}).
Thus, we treat a periodic solution of the form 
\begin{equation}
\left(r(t),\theta(t)\right)=\left(\rho,\omega t\right),\label{eq:ps}
\end{equation}
where $\omega\in\mathbb{R}$ is an angular velocity.


From (\ref{eq:DDE_PC}), the 
periodic solution satisfies\begin{subequations}\label{eq:Eqcond}
\begin{align}
1 & =\rho^{2}\cos\omega\tau,\label{eq:Eqcond1}\\
\omega & =1+\rho^{2}\sin\omega\tau.\label{eq:Eqcond2}
\end{align}
\end{subequations}
Remark that for the special case $\tau=0$, from
(\ref{eq:Eqcond}), we obtain 
\begin{equation}
\left(\rho,\omega\right)=\left(1,1\right), \label{eq:psoltau=00003D0}
\end{equation}
for periodic solution $\left(r(t),\theta(t)\right)=\left(1,t\right)$
corresponding to 
the system (\ref{eq:DDE_PC}) with $\tau=0$.

From (\ref{eq:Eqcond}) one has $\cos\omega\tau>0$. 
Then we have
the following equations\begin{subequations}\label{eq:Eqcond_a} 
\begin{align}
\omega-1 & =\tan\omega\tau,\label{eq:Eqcond_a1}\\
1+\left(\omega-1\right)^{2} & =\rho^{4}.\label{eq:Eqcond_a2}
\end{align}
\end{subequations}

\begin{figure}
\begin{center}
\subfloat[\label{fig:bif1}The roots of the equation (\ref{eq:Eqcond_a1}) are
given as intersections of two graphs $y=\omega-1$ and $y=\tan(\omega\tau)$.]{\includegraphics[scale=0.4]{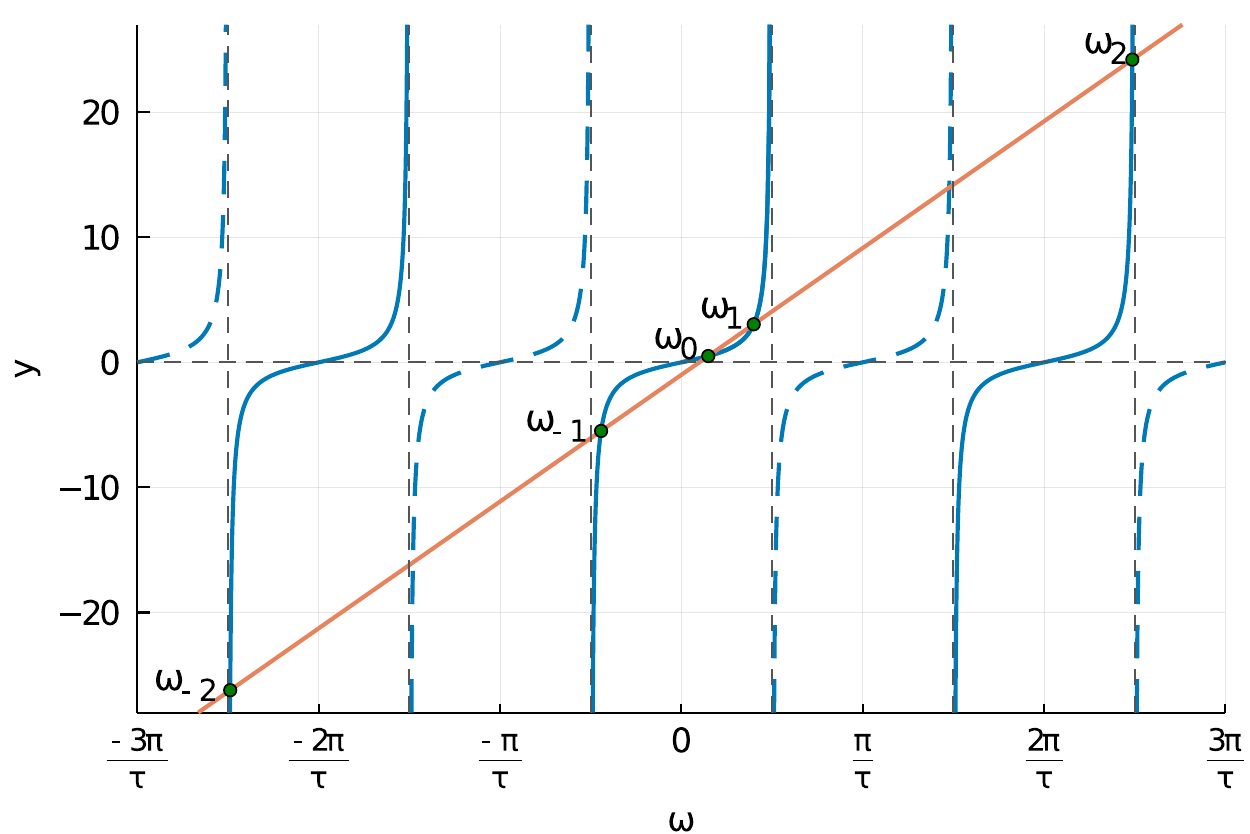}}

\subfloat[\label{fig:bif2} Branches of the periodic solutions in $\left(\tau,\omega\right)$
plane. The branch $\omega_{0}$ emerges from $\omega=1$ at $\tau=0$
and the branch $\omega_{1}$ emerges from $\omega=+\infty$ at $\tau=0$.
Those branches $\omega_{0}$ and $\omega_{1}$ meet at $\tau=\tau^{*}$.
The branches $\omega_{j}\ j\in\mathbb{Z}\setminus\left\{ 0\right\} $
appear from either $+\infty$ or $-\infty$ at $\tau=0$.]{\includegraphics[scale=0.4]{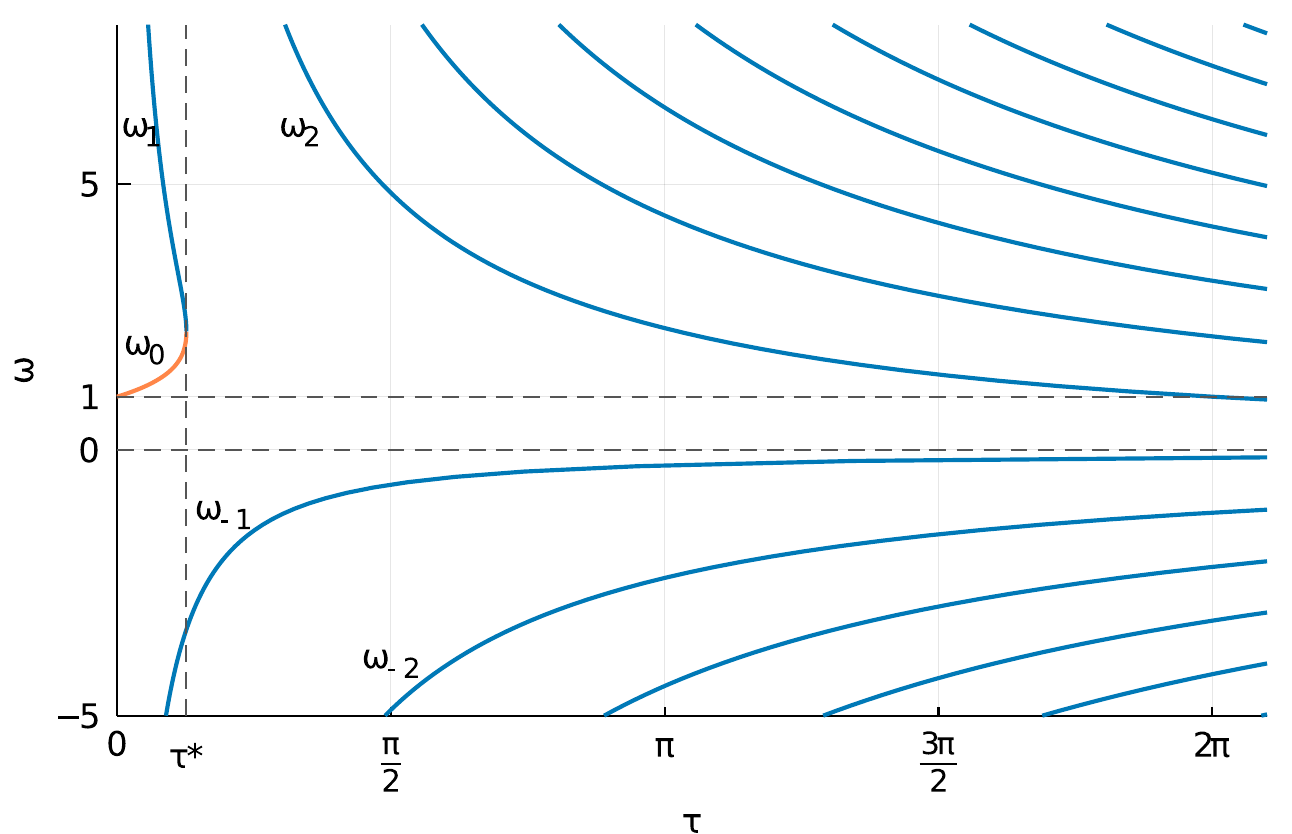}}

\caption{}\label{fig:bifall}
\end{center}
\end{figure}

In Figure \ref{fig:bifall}(a) we plot the functions $y=\omega-1$ and
$y=\tan\omega\tau$. 
For $\omega \tau \in (-\frac{\pi}{2}+2n\pi,\frac{\pi}{2}+2n\pi) (n\in \mathbb{Z}), $ 
intersections of the two functions correspond to 
the roots of (\ref{eq:Eqcond_a1}) satisfying $\cos \omega\tau>0$. 
First we study the roots of (\ref{eq:Eqcond_a1})
for $0\leq\omega\tau<\frac{\pi}{2}$. The implicit function (\ref{eq:Eqcond_a1})
for $0\leq\omega\tau<\frac{\pi}{2}$ defines a function 
\[
\overline{\tau}(\omega)=\frac{\arctan\left(\omega-1\right)}{\omega},
\]
which attains a unique maximum at $\omega=\omega^{*}$ where $\overline{\tau}^{\prime}(\omega^{*})=0$.
One can compute that 
\[
\overline{\tau}^{\prime}(\omega)=\frac{\frac{\omega}{1+\left(\omega-1\right)^{2}}-\arctan\left(\omega-1\right)}{\omega^{2}},
\]
from which we can numerically compute $\omega^{*}\approx2.22913\cdots$.
Let $\tau^{*}=\overline{\tau}(\omega^{*})$. Then we numerically obtain
\begin{equation}
\tau^{*}\approx0.398284\cdots.\label{eq:tau*}
\end{equation}

\begin{proposition}
\label{prop:ps01}For $0<\tau\leq\tau^{*}$ the equation (\ref{eq:Eqcond})
has two roots which we denote by $\omega_{0}(\tau)$ and $\omega_{1}(\tau)$
such that
\begin{itemize}
\item $\omega_{0}(\tau)\leq\omega^{*}\leq\omega_{1}(\tau),$
\item $\lim_{\tau\downarrow0}\omega_{0}(\tau)=1,\ \lim_{\tau\downarrow0}\omega_{1}(\tau)=\infty$
and
\item $\omega_{0}(\tau^{*})=\omega_{1}(\tau^{*})=\omega^{*}$.
\end{itemize}
\end{proposition}

For $\tau<\tau^{*}$ there is at least $2$ periodic solutions
\[
\left(r(t),\theta(t)\right)=\left(\rho_{j},\omega_{j}t\right),\quad j\in \left\{ 0,1\right\},
\]
where 
\begin{equation}
\rho_{j}:=\left(1+\left(\omega_{j}-1\right)^{2}\right)^{\frac{1}{4}}.
\label{eq:rhoj01}
\end{equation}

From Proposition \ref{prop:ps01} and (\ref{eq:rhoj01}), one sees
that $\lim_{\tau\downarrow0}\rho_{0}=1$ and $\lim_{\tau\downarrow0}\rho_{1}=\infty$.
Thus one periodic solution is a continuation of the periodic solution
(\ref{eq:psoltau=00003D0}) at $\tau=0$ and one periodic solution
emerges from infinity. 

For $\tau>0$, one can see that the equation (\ref{eq:Eqcond}) has
infinitely many roots. It is elementary to prove the following result,
thus we omit the proof. See also Figure \ref{fig:bifall}(a).
\begin{lemma}
For $j\in\mathbb{Z}\setminus\left\{ 0,1\right\} $, 
the equation (\ref{eq:Eqcond_a1}) with $\cos \omega\tau>0$ 
has exactly one root $\omega_{j}$ on the following each interval 
\begin{align*}
\omega_{j} \in
\left(\frac{2\pi}{\tau}\left(j-1\right),\frac{2\pi}{\tau}\left(j-1\right)
+\frac{\pi}{2\tau}\right) & ,\ j=2,3,\cdots,\\
\omega_{j} \in
\left(\frac{2\pi}{\tau}\left(j+1\right)-\frac{\pi}{2\tau},
\frac{2\pi}{\tau}\left(j+1\right)\right) & ,\ j=\cdots,-2,-1.
\end{align*}
\end{lemma}



\begin{figure}
\begin{center}
\includegraphics[scale=0.4]{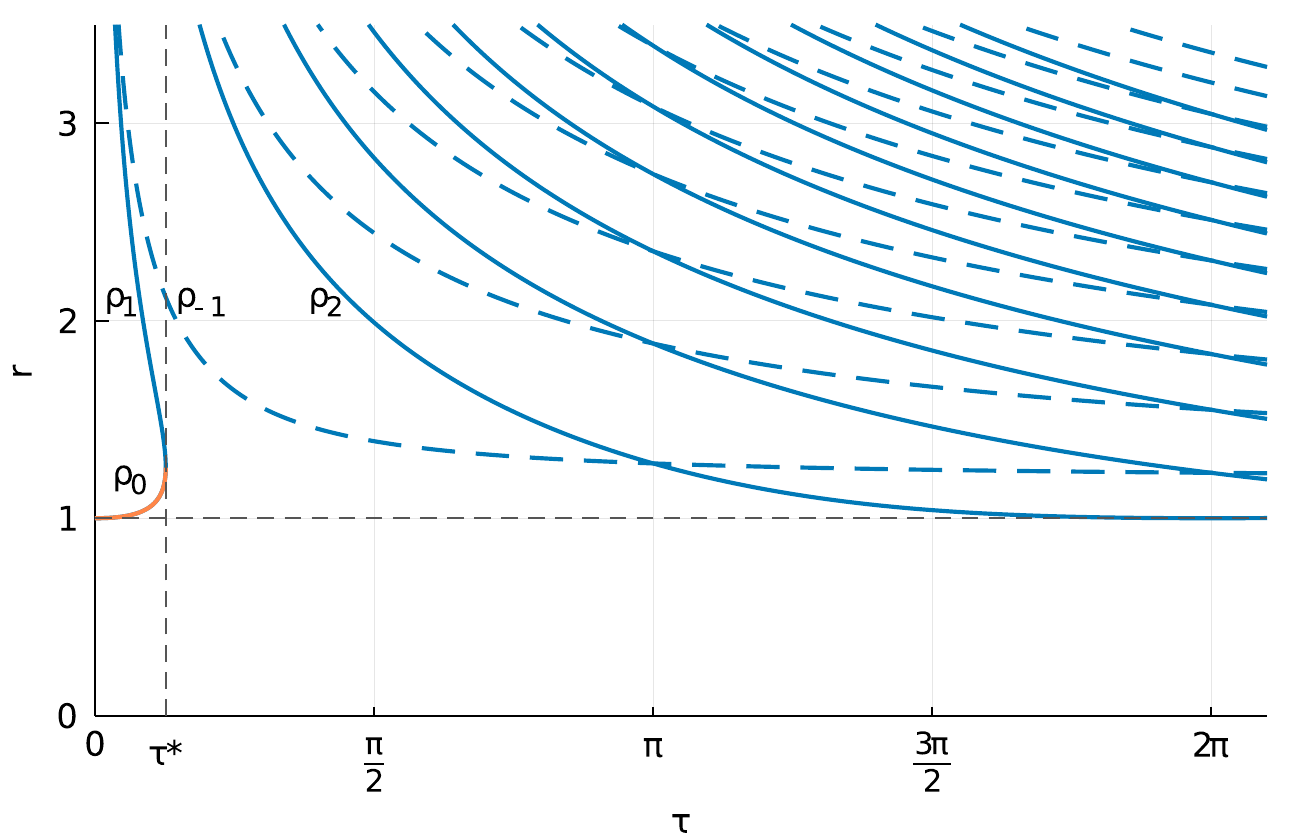}\caption{\label{fig:bifr} 
Branches of the periodic solutions in $\left(\tau,r\right)$
plane. The branch $\rho_{0}$ emerges from $r=1$ at $\tau=0$ and
the branch $\rho_{1}$ emerges from $r=+\infty$ at $\tau=0$. Those
branches $\rho_{0}$ and $\rho_{1}$ meet at $\tau=\tau^{*}$. The
branches $\rho_{j}\ j\in\mathbb{Z}\setminus\left\{ 0\right\} $ appear
from $+\infty$ at $\tau=0$ and exist for $\tau>0$. Dashed curve
denotes the periodic solutions with $\omega_{j}<0$ (clockwise periodic
solutions).}
\end{center}
\end{figure}

Therefore, we obtain roots $\omega_{j}$ for $j\in\mathbb{Z}$ if
$\tau\leq\tau^{*}$ 
and for $j\in\mathbb{Z}\setminus\left\{ 0,1\right\} $ if $\tau>\tau^{*}$
of the equation (\ref{eq:Eqcond_a1}). In Figure \ref{fig:bifall}(b) we plot
the branches for $\omega$ as a function of $\tau$. Once $\omega$
of the periodic solution (\ref{eq:ps}) is given, $\rho$ is determined
from (\ref{eq:Eqcond_a2}). The radius $\rho$ can be determined
as $\rho_{j}$ given as, similar to (\ref{eq:rhoj01}), 
\[
\rho_{j}=\left(1+(\omega_{j}-1)^{2}\right)^{\frac{1}{4}}\geq1.
\]
This implies that the periodic solutions of the form (\ref{eq:ps})
has larger radius than the unit circle which is the trajectory of
the periodic solution for $\tau=0$. We also note that there is a
root $\rho=1$ if $\tau=2n\pi,\ n\in\mathbb{N}$. In Figure \ref{fig:bifr}
we plot the branches for the radius as a function of $\tau$.

Summarizing the above findings, we obtain the following result for
the existence of the periodic solution (\ref{eq:ps}) for the system
(\ref{eq:DDE_PC}). The result is not intuitive and not expected that
the delay induces many periodic solutions from the non-delay system
(\ref{eq:DDE}).
\begin{theorem}\label{thm:existence-periodicsol}
For each $\tau>0$, the equation (\ref{eq:Eqcond}) has infinitely
many roots, which are countable. Thus the system (\ref{eq:DDE_PC})
has infinitely many periodic solutions, which are countable, of the
form (\ref{eq:ps}).
\end{theorem}

\section{Analysis of the characteristic equation for the stability of the periodic
solutions}

To analyze stability of the periodic solution obtained in Section
3, we study a system of a delay differential equation and an integral
equation, employing the principle of linearized stability for the
coupled systems of renewal equations and delay differential equations
\cite{Diekmann:2008}.

Let $v(t)=\theta^{\prime}(t)$. From the system (\ref{eq:DDE_PC}),
we get the following system of a delay differential equation and a
renewal equation\begin{subequations}\label{eq:DE}
\begin{align}
r^{\prime}(t) & =r(t)\left(1-r(t)r(t-\tau)\cos\int_{t-\tau}^{t}v(s)ds\right),\label{eq:DE1}\\
v(t) & =1+r(t)r(t-\tau)\sin\int_{t-\tau}^{t}v(s)ds.\label{eq:DE2}
\end{align}
\end{subequations}The system (\ref{eq:DE}) has equilibria
\begin{equation}
\left(r(t),v(t)\right)=\left(\rho_{j},\omega_{j}\right),\ j\in\mathbb{N}\label{eq:equilibrium}
\end{equation}
corresponding to the periodic solutions 
\[
\left(r(t),\theta(t)\right)=\left(\rho_{j},\omega_{j}t\right),\ j\in\mathbb{N}
\]
of the system (\ref{eq:DDE_PC}), given as in Section 3.

In the Appendix A, by linearization of the system (\ref{eq:DE}) at
the equilibrium (\ref{eq:equilibrium}), we obtain the following characteristic
equation
\begin{equation}
0=\lambda+2e^{-\lambda\tau}-\rho^{4}\int_{0}^{2\tau}e^{-\lambda s}ds,\ \lambda\in\mathbb{C},\label{eq:cheq}
\end{equation}
with $\rho=\rho_{j},\ j\in\mathbb{N}$. Note that we have a family
of the characteristic equations (\ref{eq:cheq}) which are indexed
by the equilibrium about which we linearize. 

Define $\eta:=\lambda\tau$. Then, from (\ref{eq:cheq}), we obtain
the following equation
\begin{equation}
0=f(\eta),\ \eta\in\mathbb{C},\label{eq:cheq_eta}
\end{equation}
where $f:\mathbb{C}\to\mathbb{C}$ defined as
\begin{equation}
f(\eta):=\eta+2\tau e^{-\eta}-\rho^{4}\tau^{2}\int_{0}^{2}e^{-\eta s}ds,\ \eta\in\mathbb{C}.\label{eq:lhch}
\end{equation}
We study (\ref{eq:cheq_eta}) to analyze the distribution of the
roots of (\ref{eq:cheq}) with respect to the imaginary axis. First
let us study the existence of real roots for the characteristic equation. 
\begin{lemma}
\label{lem:crit}For each $j\in\mathbb{N}$, the characteristic equation
(\ref{eq:cheq_eta}) has 
\begin{enumerate}
\item a unique negative real root if $\rho_{j}^{4}\tau<1$,
\item a root $\lambda=0$ if $\rho_{j}^{4}\tau=1$, and 
\item a unique positive real root if $\rho_{j}^{4}\tau>1$.
\end{enumerate}
\end{lemma}

\begin{proof}
Consider the function $f(\eta)$ for $\eta\in\mathbb{R}$. Note that
$f(0)=2\tau\left(1-\rho^{4}\tau\right).$ Since it follows that 
\[
\lim_{\eta\to-\infty}e^{\eta}\int_{0}^{2}e^{-\eta s}ds=\lim_{\eta\to-\infty}\frac{e^{\eta}-e^{-\eta}}{\eta}=\infty,
\]
one obtains
\[
\lim_{\eta\to\infty}f(\eta)=\infty,\ \lim_{\eta\to-\infty}f(\eta)=-\infty.
\]
Since $f$ is a continuous function, if $\rho_{j}^{4}\tau<1$ then
there is a negative real root and if $\rho_{j}^{4}\tau>1$ then there
is a positive real root. We also see that if $\rho_{j}^{4}\tau=1$
then there exists a root $0$ for the function $f$.

For the uniqueness of the root, we study the equation $0=\eta f(\eta)$.
Here 
\begin{align*}
\eta f(\eta) & =\eta^{2}+2\tau e^{-\eta}\eta-\rho^{4}\tau^{2}\left(1-e^{-2\eta}\right)\\
 & =\left(\eta+\tau e^{-\eta}\right)^{2}-\tau^{2}e^{-2\eta}-\rho^{4}\tau^{2}\left(1-e^{-2\eta}\right).
\end{align*}
Let us define 
\begin{align*}
g_{1}(\eta) & =\left(\eta+\tau e^{-\eta}\right)^{2},\\
g_{2}(\eta) & =\tau^{2}e^{-2\eta}+\rho^{4}\tau^{2}\left(1-e^{-2\eta}\right).
\end{align*}
We consider an intersection of $g_{1}$ and $g_{2}$. It is easy to
see that $g_{1}(0)=g_{2}(0)=\tau^{2}$. We compute 
\begin{align*}
g_{1}^{\prime}(\eta) & =2\left(\eta+\tau e^{-\eta}\right)\left(1-\tau e^{-\eta}\right),\\
g_{2}^{\prime}(\eta) & =2\tau^{2}\left(\rho^{4}-1\right)e^{-2\eta}>0,\\
g_{2}^{\prime\prime}(\eta) & =-4\tau^{2}\left(\rho^{4}-1\right)<0.
\end{align*}
Therefore, one sees that $g_{1}$ is a downward-convex function (attaining
minimum at $\eta=\log\tau$) and $g_{2}$ is upward-convex function.
Hence, the intersection of the functions $g_{1}$ and $g_{2}$ except
$0$ is unique. Thus we obtain the conclusion.
\qed
\end{proof}
For each $j\in\mathbb{N}$ we obtain the estimation of $\rho_{j}^{4}\tau$
with respect to $1$ as follows. 
\begin{lemma}
\label{lem:whichbranch}The following statements are true.
\begin{itemize}
\item It holds $\rho_{0}^{4}\tau<1$ and $\rho_{1}^{4}\tau>1$ for $\tau<\tau^{*}$
and $\rho_{0}^{4}\tau=\rho_{1}^{4}\tau=1$ for $\tau=\tau^{*}$, and 
\item it holds $\rho_{j}^{4}\tau>1\ \left(j\in\mathbb{N}\setminus\left\{ 0,1\right\} \right)$
for $\tau>0$. 
\end{itemize}
\end{lemma}

\begin{proof}
For any $j\in\mathbb{Z}$, from the condition (\ref{eq:Eqcond_a}),
\[
\rho_{j}^{4}=1+\tan^{2}\omega_{j}\tau=\frac{1}{\cos^{2}\omega_{j}\tau}.
\]
Thus 
\[
\rho_{j}^{4}\tau=\frac{d}{d\omega}\tan(\omega\tau)|_{\omega=\omega_{j}}
\]
holds. One can see that 
\[
\frac{d}{d\omega}\tan(\omega\tau)|_{\omega=\omega_{0}}<1<\frac{d}{d\omega}\tan(\omega\tau)|_{\omega=\omega_{1}}
\]
for $\tau<\tau^{*}$ and 
\[
\frac{d}{d\omega}\tan(\omega\tau)|_{\omega=\omega_{0}}=1=\frac{d}{d\omega}\tan(\omega\tau)|_{\omega=\omega_{1}}
\]
for $\tau=\tau^{*}$. Thus we obtain the first statement. For $j\in\mathbb{Z}\setminus\left\{ 0,1\right\} $
it is clear that 
\[
\rho_{j}^{4}\tau=\frac{d}{d\omega}\tan(\omega\tau)|_{\omega=\omega_{j}}>1.
\]
Thus we obtain the conclusion.
\qed
\end{proof}
From Lemmas \ref{lem:crit} and \ref{lem:whichbranch} 
and the principle
of linearized stability \cite{Diekmann:2008}, we obtain the following
result concerning instability of the periodic solutions.
\begin{theorem}
The periodic solution of the form (\ref{eq:ps}) for $j=1$ is unstable
for $\tau<\tau^{*}$. The periodic solution of the form (\ref{eq:ps})
for $j\in\mathbb{Z}\setminus\left\{ 0,1\right\} $ is unstable. 
\end{theorem}

Let us consider stability of the periodic solution of the form (\ref{eq:ps})
for $j=0$. For $j=0$ the characteristic equation (\ref{eq:cheq_eta})
may have an imaginary root with positive real part. Thus so far we
cannot determine stability of the periodic solution from Lemma \ref{lem:crit}
and \ref{lem:whichbranch}. 
Below we exclude this possibility to conclude
that the periodic solution of the form (\ref{eq:ps}) for $j=0$ is
asymptotically stable for $\tau<\tau^{*}$.

First let us show that there is a compact region in the complex plane
for the existence of a root with positive real part.
\begin{lemma}
\label{lem:smalltau}Let $j=0$. Suppose that the characteristic equation
(\ref{eq:cheq_eta}) has a root $\eta$ with $\mathrm{Re}\,\eta>0$.
Then every root $\eta$ with $\mathrm{Re}\,\eta>0$ satisfies
\begin{equation}
\left|\eta\right|\le2\tau+2\rho_{0}^{4}\tau^{2}\frac{1}{\left|\eta\right|}.\label{eq:est}
\end{equation}
Thus for any $\varepsilon>0$ there exists $\delta>0$ such that $\tau<\delta$
implies $\left|\eta\right|<\varepsilon$.
\end{lemma}

\begin{proof}
Consider the equation (\ref{eq:cheq_eta}) for $\eta$ such that $\mathrm{Re}\,\eta>0$.
We have 
\[
\int_{0}^{2}e^{-\eta s}ds=\frac{1-e^{-2\eta}}{\eta}.
\]
Thus from the characteristic equation (\ref{eq:cheq_eta}), we obtain
\[
\eta=-2\tau e^{-\eta}+\rho_{0}^{4}\tau^{2}\frac{1-e^{-2\eta}}{\eta}.
\]
By the straightforward estimation $\left|e^{-\eta}\right|<1$ and
$\left|1-e^{-2\eta}\right|\leq 1+ \left|e^{-2\eta}\right|<2$ using
  $\mathrm{Re}\:\eta>0$, we obtain
the estimation (\ref{eq:est}).
From the inequality (\ref{eq:est}), one gets 
\[
\left|\eta\right|\leq \tau+\sqrt{2\tau+\tau^2},
\]
from which we obtain the conclusion.
\qed

\end{proof}
Substituting $\eta=\mu+i\nu,\ \left(\eta,\nu\right)\in\mathbb{R}^{2}$
into (\ref{eq:cheq_eta}), we obtain the following two equations.
\begin{subequations}\label{eq:CHEQ_RI}
\begin{align}
0 & =\mu+2\tau e^{-\mu}\cos\nu-\rho_{j}^{4}\tau^{2}\int_{0}^{2}e^{-\mu s}\cos\nu s\, ds,\label{eq:CHEQ_R}\\
0 & =\nu-2\tau e^{-\mu}\sin\nu+\rho_{j}^{4}\tau^{2}\int_{0}^{2}e^{-\mu s}\sin\nu s\, ds.\label{eq:CHEQ_I}
\end{align}
\end{subequations}From (\ref{eq:CHEQ_I}) it is immediate to obtain
the following lemma. We omit the proof.
\begin{lemma}
\label{lem:conj}
If $\eta=\mu+i\nu$ is a root of the characteristic equation (\ref{eq:cheq_eta}),
then so is its conjugate $\overline{\eta}=\mu-i\nu$.
\end{lemma}

From Lemma \ref{lem:conj}, it is sufficient to consider the root $\eta=\mu + i\nu$ of the characteristic equation (\ref{eq:cheq_eta}) with $\nu>0$. Our next aim is to show the following result.
\begin{lemma}
Let $j=0$. There exists $\varepsilon$ such that if $\tau<\varepsilon$
then there is no root with positive real part for the characteristic
equation (\ref{eq:cheq_eta}). 
\end{lemma}

\begin{proof}
Consider the branch $j=0$ for $\tau<\tau^{*}$. Suppose that there
exists a root $\eta=\mu+i\nu$ with $\mu>0$. From Lemma \ref{lem:smalltau}, there exists $\tau$ such that $\sin\nu s\geq 0$ for $s\in[0,2]$. Thus, for this $\tau$, one has $\int_{0}^{2}e^{-\mu s}\sin\nu sds\geq0$.
Furthermore, for sufficiently small $\tau$, 
\[
\nu-2\tau e^{-\mu}\sin\nu=\nu\left(1-2\tau e^{-\mu}\frac{\sin\nu}{\nu}\right)>0
\]
for $\nu>0$. Therefore, we obtain a contradiction to (\ref{eq:CHEQ_RI}).
Thus we obtain the conclusion.
\qed
\end{proof}
\begin{lemma}
Let $j=0$. The characteristic equation (\ref{eq:cheq_eta}) does
not have a purely imaginary root.
\end{lemma}

\begin{proof}
Assume that there exists an imaginary root $\eta=i\nu,\ \nu>0$ for
the characteristic equation (\ref{eq:cheq_eta}). Then, from (\ref{eq:CHEQ_I}),
it holds that \begin{subequations}\label{eq:CHEQ_RI2}
\begin{align}
0 & =2\tau\cos\nu-\rho_{0}^{4}\tau^{2}\int_{0}^{2}\cos\nu s\, ds,\label{eq:CHEQ_R2}\\
0 & =\nu-2\tau\sin\nu+\rho_{0}^{4}\tau^{2}\int_{0}^{2}\sin\nu s\, ds.\label{eq:CHEQ_I2}
\end{align}
\end{subequations}One can see that 
\begin{align*}
\int_{0}^{2}\cos\nu s\, ds & =\frac{\sin2\nu}{\nu}=\frac{2\sin\nu\cos\nu}{\nu},\\
\int_{0}^{2}\sin\nu s\, ds & =\frac{1-\cos2\nu}{\nu}=\frac{2\sin^{2}\nu}{\nu}.
\end{align*}
Then, from (\ref{eq:CHEQ_RI2}), we obtain \begin{subequations}\label{eq:CHEQ_RI3}
\begin{align}
0 & =\cos\nu\left(1-\rho_{0}^{4}\tau\frac{\sin\nu}{\nu}\right),\label{eq:CHEQ_R2-1}\\
0 & =\nu\left(1-2\tau\frac{\sin\nu}{\nu}+2\rho_{0}^{4}\tau^{2}\left(\frac{\sin\nu}{\nu}\right)^{2}\right).\label{eq:CHEQ_I2-1}
\end{align}
\end{subequations}Note that we have $\rho_{0}\geq1$ from (\ref{eq:rhoj01}).
Since, for any $a\in\mathbb{R}$, it holds that $1-2\tau a+2\tau^{2}a^{2}>0,$
we see that 
\[
1-2\tau\frac{\sin\nu}{\nu}+2\rho_0^{4}\tau^{2}\left(\frac{\sin\nu}{\nu}\right)^{2}\geq1-2\tau\frac{\sin\nu}{\nu}+2\tau^{2}\left(\frac{\sin\nu}{\nu}\right)^{2}>0,
\]
which is a contradiction to (\ref{eq:CHEQ_I2-1}). Therefore, we obtain
the conclusion.
\qed
\end{proof}
Therefore, from the application of the Rouche's theorem (see e.g.
Lemma 2.8 in Chapter XI of \cite{Diekmann:1995}), we obtain the following
conclusion.
\begin{theorem}
\label{thm:psolst}If $\tau<\tau^{*}$, then the periodic solution
of the form (\ref{eq:ps}) with $j=0$ is asymptotically stable and
the periodic solution $j\not=0$ is unstable. If $\tau>\tau^{*}$
then every periodic solution (\ref{eq:ps}) is unstable.
\end{theorem}

From Theorems \ref{thm:existence-periodicsol} and \ref{thm:psolst}, 
we complete the proof of Theorem \ref{thm:periodic}.

\section{Numerical simulations}  \label{Numerics}

In this section we demonstrate numerical solutions of the delay differential
equation (\ref{eq:DDE}) with the special initial conditions $x(t)=y(t)=\delta,\ t\in\left[-\tau,0\right]$, where
$\delta\in\mathbb{R}$. 

\begin{figure}
\begin{center}
\includegraphics[scale=0.4]{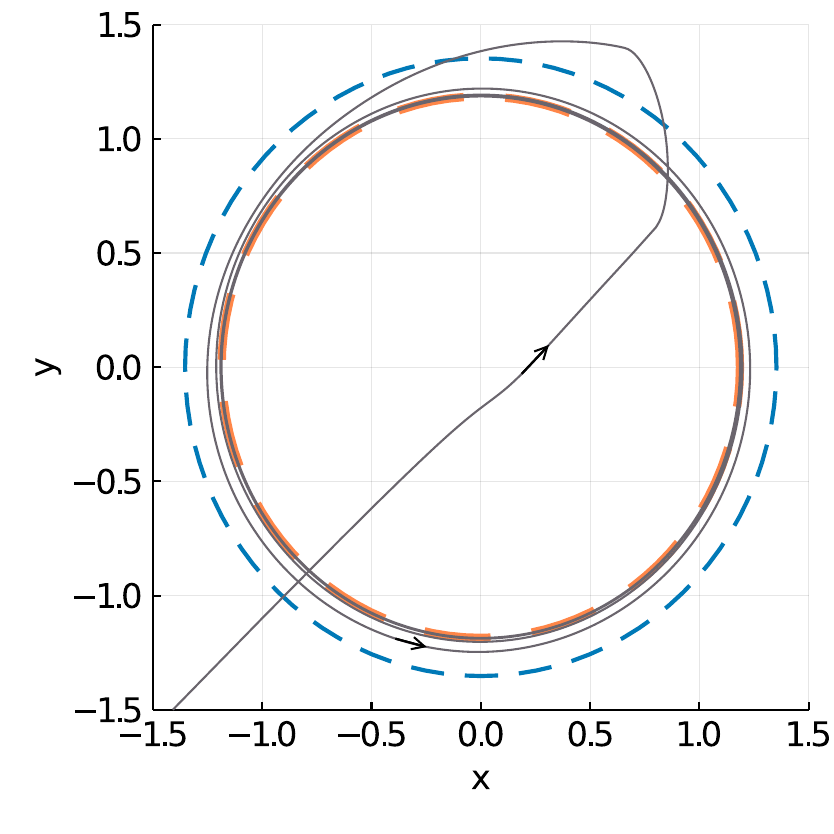}
\caption{\label{fig:stab}
Illustration of a trajectory of the solution with
$\delta=-36$ in $(x,y)$-plane for $\tau=0.392$. 
The asymptotic
stable solution attracts the solution with the initial condition $\delta=-36$.
The trajectory of the asymptotically stable periodic solution and
of the unstable periodic solution are illustrated as the dashed orange
circle and as the dashed blue circle, respectively.}
\end{center}
\end{figure}

\begin{figure}
\begin{center}
\includegraphics[scale=0.4]{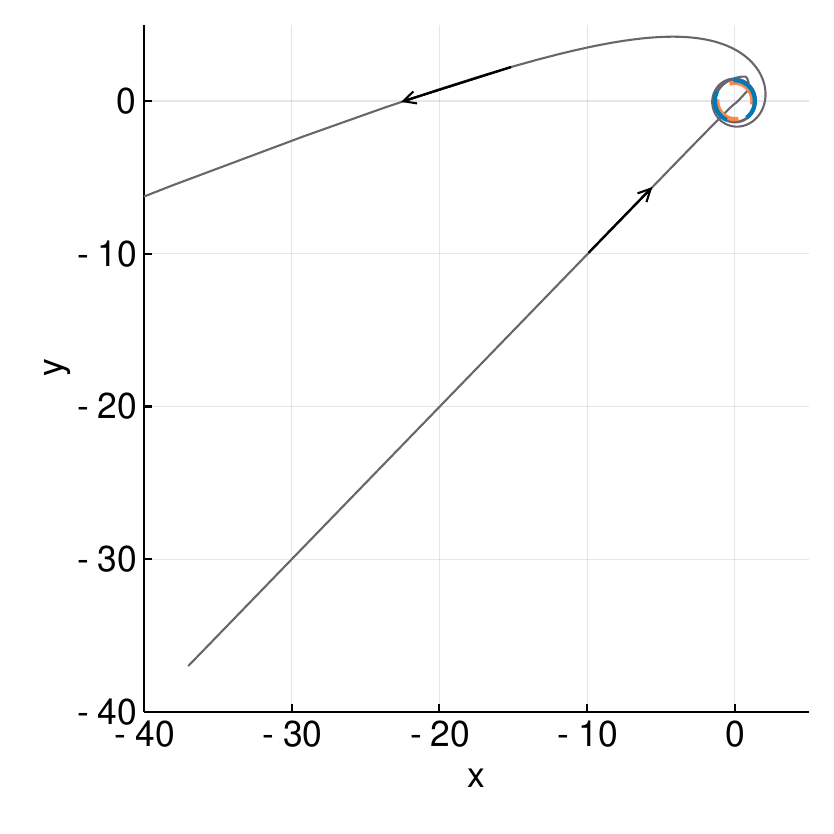}\includegraphics[scale=0.4]{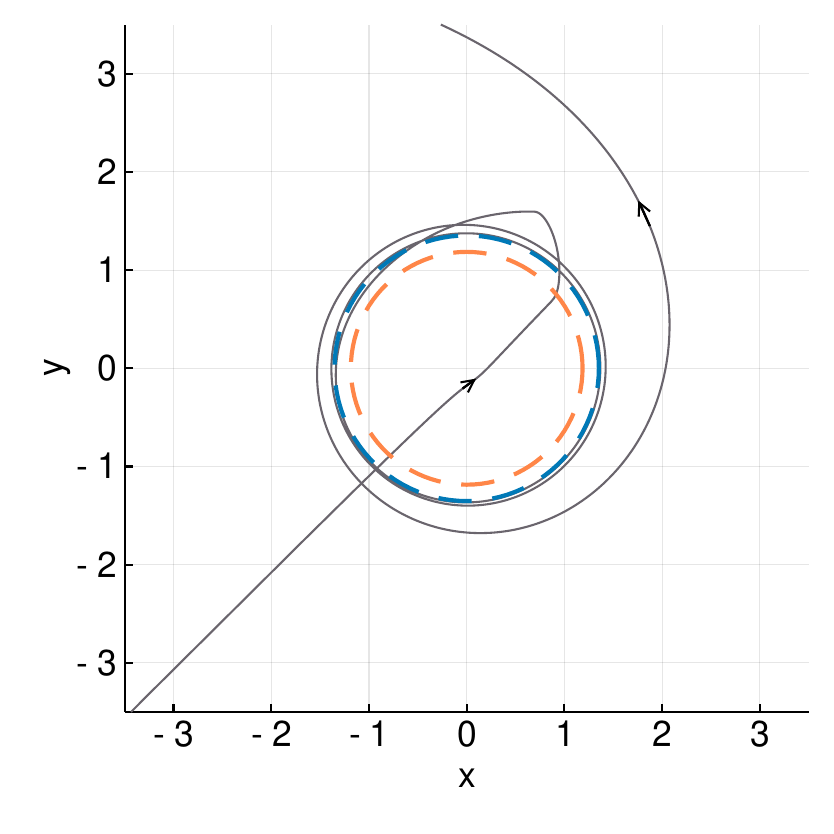}
\caption{\label{fig:unst}
Illustration of a trajectory of the solution with
$\delta=-37$ in $(x,y)$-plane for $\tau=0.392$. 
(The right figure is close-up view of the left figure near the origin.)
The solution blows up,
after winding around the unstable periodic solution (dashed blue circle).
The trajectory of the asymptotically stable periodic solution and
of the unstable periodic solution are illustrated as the dashed orange
circle and as the dashed blue circle, respectively.}
\end{center}
\end{figure}

First we fix $\tau=0.392$. For $\tau<\tau^{*}\approx0.398284\cdots$,
we show that the periodic solution of the form (\ref{eq:ps}) with
$j=0$ is asymptotically stable (Theorem \ref{thm:psolst}). From
(\ref{eq:Eqcond_a}), we can compute the radius of the stable periodic
solution as $r=\rho_{0}\approx1.18547\cdots$. In Figure \ref{fig:stab},
a trajectory of the solution for $\delta=-36$ is plotted in $(x,y)$
plane, which shows that the periodic solution with the radius $r=\rho_{0}\approx1.18547\cdots$
attracts the solution. Observe that there is an unstable periodic
solution in the vicinity of the stable periodic solution (in Figure
\ref{fig:stab}, the trajectory of the asymptotically stable periodic
solution and of the unstable periodic solution are illustrated as
the dashed orange circle and as the dashed blue circle, respectively).
In Figure \ref{fig:unst}, a trajectory of the solution with the initial
condition $\delta=-37$ is plotted. In this case, the solution winds
around the unstable periodic solution (dashed blue circle) and then
leaves and goes far away. 

\begin{figure}
\begin{center}
\includegraphics[scale=0.5]{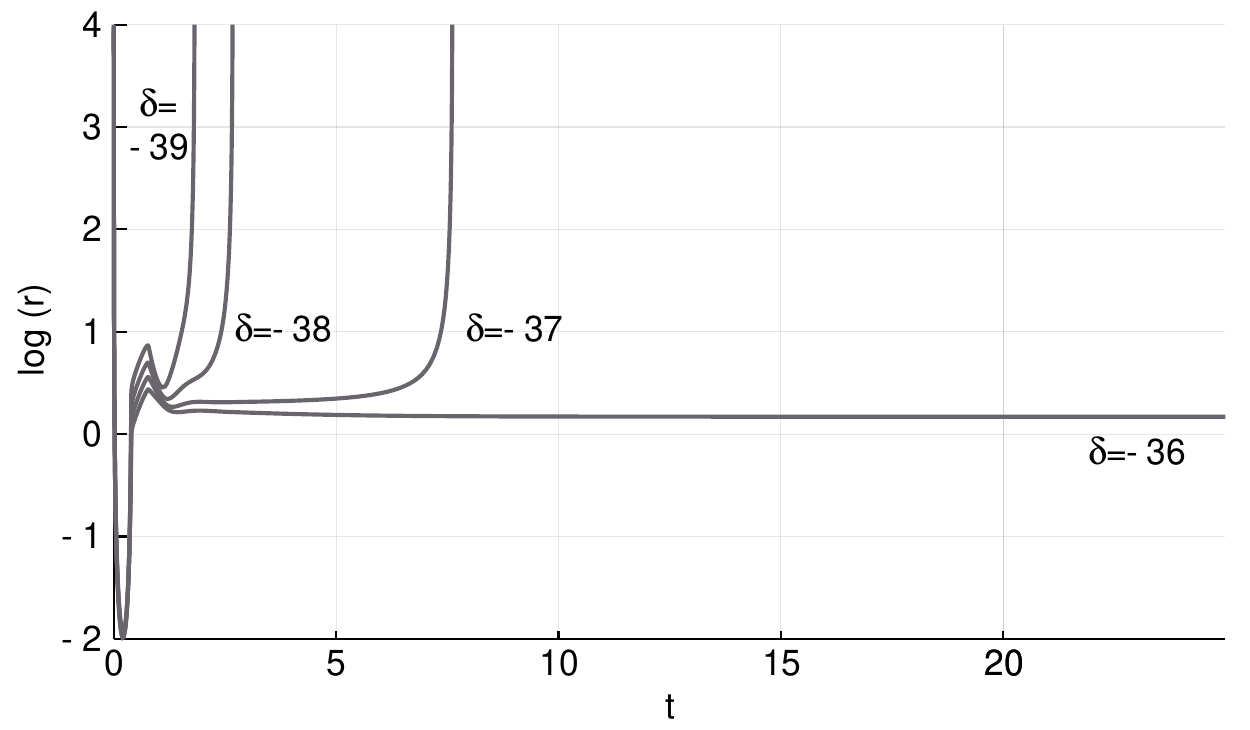}\caption{\label{fig:growth}We plot the growth of $\log r(t)$ for several
$\delta$. Here $\tau$ is fixed ($\tau=0.392)$. The numerical simulation
suggests that larger $\left|\delta\right|$ causes the solution to blow up faster.}
\end{center}
\end{figure}

In Figure \ref{fig:growth}, we plot $\log r(t)$ for several initial
conditions $\delta$. The numerical experiment suggests that the solution
exists globally for $-36<\delta\leq0$ and blows up for $\delta<-37$:
in short, the solution blows up for large $\left|\delta\right|$.
We can numerically observe many blow-up solutions, which blow up 
even after $t=\tau/2$. In this paper, we prove the existence of blow-up
solutions, which blow up in the time interval $(0,\tau/2)$. The numerical
simulation suggest that many solutions blow up in finite times. 

\begin{figure}
\begin{center}
\includegraphics[scale=0.5]{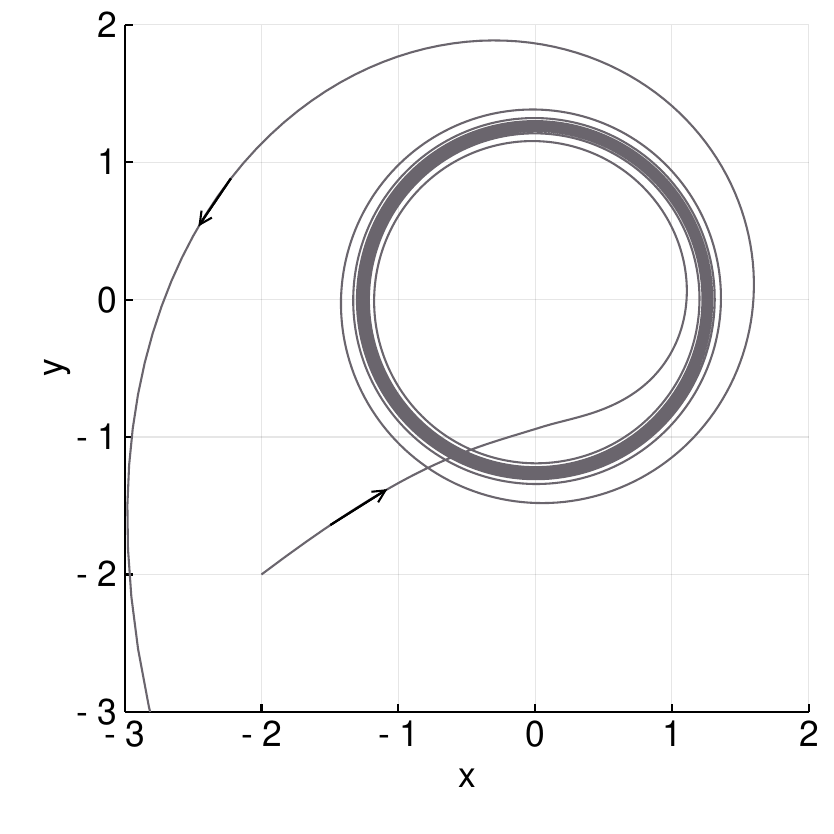}\caption{\label{fig:trans}Transient behavior of the solution for $\tau=0.3985>\tau^{*}$.
The solution stays around the origin for a while, then it blows up 
in a finite time.}
\end{center}
\end{figure}

We also numerically compute the solution for $\tau>\tau^{*}$. Figure
\ref{fig:trans} shows a transient behavior of a solution. The solution
stays around the origin for a while and then it blows up in a finite
time. Although for $\tau=0.3985>\tau^{*}$, there does not exist periodic
solution for $j=0$ and $j=1$ around the origin, those periodic solutions
may indirectly affect the transient behavior of the solution.

\section{Discussion}

In this paper we show an example of a blow-up phenomenon in a planar
system of delay differential equations. In our system, the delay
completely changes the system: many blow-up solutions and periodic
solutions suddenly appear, no matter how small the length of the delay
is. In neutral delay differential equations, it is known that arbitrary
small delay can destabilize the system (see Chapter 1.7 in \cite{Hale:1993}).
Here we find that arbitrary small delay can induce blow-up solutions 
in delay differential equations\@. Numerical simulations suggest
that many solutions either tend to the stable periodic solution or
blow up in a finite time. It is not obvious if more complicated
solution behavior exists. The transient behavior observed in Figure
\ref{fig:trans} is interesting, as it looks that the solution tries
to find a stable periodic solution which does not exist in this parameter
setting.

Many results concerning the existence of the blow-up solutions are
available in Volterra integral equations (see \cite{Bandle:1998,Brunner:2012}
and references therein). 
When delay differential equations can be
formulated as Volterra integral equations, we may apply the blow-up
results for Volterra integral equations to delay differential equations,
see e.g. \cite{Appleby:2017,Yang:2013} for a relation between Volterra type
integral equations and integro-differential equations. 
However, 
we cannot apply the results for Volterra integral equations to 
a system of integral equations which is rewritten formally from
 our target system of delay differential equations.
Since delay differential equations form an infinite dimensional dynamical system,
it is also not straightforward to apply the results established in
ordinary differential equations. 

In the context of mathematical modelling, our example suggests that
arbitrary small delay can be responsible for a drastic change of the
dynamics, thus one should be careful when ignoring small delay. Other
blow-up mechanisms in delay differential equations will be explored
in our future work.

\section*{Acknowledgements}

The authors thank the referees for careful reading and invaluable comments
on our manuscript.

\appendix

\section{The characteristic equation}

Here, by linearization of the system (\ref{eq:DE}) about the equilibrium
(\ref{eq:equilibrium}), we derive the characteristic equation (\ref{eq:cheq}),
which characterizes stability of the periodic solutions. Fixing $j\in\mathbb{Z}$,
we let 
\[
s_{j}(t)=r(t)-\rho_{j},\ u(t)=v(t)-\omega_{j}.
\]
Then, applying the Taylor expansion, we get
\begin{align*}
\cos\int_{t-\tau}^{t}v(s)ds & =\cos\int_{t-\tau}^{t}\left(u(s)+\omega_{j}\right)ds=\cos\omega_{j}\tau-\sin\omega_{j}\tau\int_{t-\tau}^{t}u(s)ds+\cdots,\\
\sin\int_{t-\tau}^{t}v(s)ds & =\sin\int_{t-\tau}^{t}\left(u(s)+\omega_{j}\right)ds=\sin\omega_{j}\tau+\cos\omega_{j}\tau\int_{t-\tau}^{t}u(s)ds-\cdots.
\end{align*}
Ignoring the higher order terms, we obtain the following linearized
equation\begin{subequations}\label{eq:LE}
\begin{align}
s^{\prime}(t) & =-s(t)-s(t-\tau)+\rho_{j}^{3}\sin\omega_{j}\tau\int_{0}^{\tau}u(t-s)ds,\label{eq:LE1}\\
u(t) & =
\rho_j\sin\omega_{j}\tau\left(s(t)+s(t-\tau)\right)+\int_{0}^{\tau}u(t-s)ds.\label{eq:LE2}
\end{align}
\end{subequations}Substituing the exponential solution $\left(s(t),u(t)\right)=e^{\lambda t}\mathbf{c}$,
where $\lambda\in\mathbb{C}$ and $\mathbf{c}\in\mathbb{C}^{2}$,
we get the following characteristic equation 
\begin{equation}
\text{det}\left(\left[\begin{array}{cc}
-1-e^{-\lambda\tau} & \rho_{j}^{3}\sin\omega\tau\int_{0}^{\tau}e^{-\lambda s}ds\\
\rho_{j}\sin\omega_{j}\tau\left(1+e^{-\lambda\tau}\right) & \int_{0}^{\tau}e^{-\lambda s}ds
\end{array}\right]-\left[\begin{array}{cc}
\lambda & 0\\
0 & 1
\end{array}\right]\right)=0.\label{eq:CHEQ0}
\end{equation}
Using (\ref{eq:Eqcond1}), we have
\[
\rho_{j}^{4}\sin^{2}\omega\tau=\rho_{j}^{4}-\rho_{j}^{4}\cos^{2}\omega\tau=\rho_{j}^{4}-1.
\]
Therefore, the characteristic equation (\ref{eq:CHEQ0}) becomes
\begin{equation}
\left(1+e^{-\lambda\tau}+\lambda\right)\left(1-\int_{0}^{\tau}e^{-\lambda s}ds\right)-\left(\rho_{j}^{4}-1\right)\left(1+e^{-\lambda\tau}\right)\int_{0}^{\tau}e^{-\lambda s}ds=0.\label{eq:cheqcomp1}
\end{equation}
The equation (\ref{eq:cheqcomp1}) can be written as 
\[
\left(1+e^{-\lambda\tau}+\lambda\right)-\int_{0}^{\tau}e^{-\lambda s}ds\left(\lambda+\rho_{j}^{4}\left(1+e^{-\lambda\tau}\right)\right)=0.
\]
Therefore, we obtain the characteristic equation (\ref{eq:cheq})
in the main text.

%
%


\begin{thebibliography}{10}
\bibitem{Appleby:2017}J.A.D. Appleby, D.D. Patterson, 
Blow-up and superexponential growth in superlinear Volterra equations. 
Disc. Cont. Dyn. Sys., 38 (8) (2017) pp. 3993--4017.

\bibitem{Bandle:1998}C. Bandle, H. Brunner, 
Blowup in diffusion equations: a survey. 
J. Comp. Appl. Math., 97 (1-2) (1998) pp. 3--22.

\bibitem{Brunner:2012}H. Brunner, Z.W. Yang, 
Blow-up behavior of Hammerstein-type Volterra integral equations. 
J. Int. Equ. Appl., 24 (4) (2012) pp. 487--512.

\bibitem{Diekmann:2008}O. Diekmann, Ph. Getto, M. Gyllenberg, 
Stability and bifurcation analysis of Volterra functional equations in the light of suns and stars. 
SIAM J. Math. Anal., 39 (4) (2008) pp. 1023--1069.

\bibitem{Diekmann:1995}O. Diekmann, S.A. van Gils, S.M. Verduyn Lunel, H.O. Walther, 
Delay Equations: Functional-, Complex-, and Nonlinear Analysis, 
Applied Mathematical Sciences (Vol. 110), Springer (1995)

\bibitem{Elias:2006}U. Elias, H. Gingold, 
Critical points at infinity and blow up of solutions of autonomous polynomial differential systems
via compactification. 
J. Math. Anal. Appl., 318 (1) (2006) pp. 305--322.

\bibitem{Erneux:2009}T. Erneux, 
Applied Delay Differential Equations. 
Springer (2009)

\bibitem{EzzinbiJazar2006}K. Ezzinbi, M. Jazar, 
Blow-up results for some nonlinear delay differential equations. 
Positivity, 10 (2006) pp. 329--341.

\bibitem{Hale:1993}J.K. Hale, S.M. Verduyn Lunel, 
Introduction to Functional Differential Equations. 
Springer, New York (1993)

\bibitem{Hirata:1999}D. Hirata, 
Blow-up for a class of semilinear integro-differential equations of parabolic type. 
Math. Meth. Appli. Sci., 22 (13) (1999) pp. 1087--1100.

\bibitem{Hirota:2006}C. Hirota, K. Ozawa, 
Numerical method of estimating the blow-up time and rate of the solution of ordinary differential 
equations - An application to the blow-up problems of partial differential equations. 
J. Comp. Appl. Math., 193 (2) (2006) pp. 614--637.

\bibitem{Kuznetsov:1998}
Y. A. Kuznetsov, 
Elements of Applied Bifurcation Theory, Springer (1998)

\bibitem{Ma:2011}J. Ma, 
Blow-up solutions of nonlinear Volterra integro-differential equations. 
Math. Comp. Model., 54 (11--12) (2011) pp. 2551--2559.

\bibitem{Matsue:2018}K. Matsue, 
On blow-up solutions of differential equations with Poincar\'e-type compactifications. 
SIAM J. Appl. Dyn. Syst., 17 (3) (2018) pp. 2249--2288.

\bibitem{Mydlarczyk}W. Mydlarczyk, 
A condition for finite blow-up time for a Volterra integral equation. 
J. Math. Anal. Appl., 181 (1) (1994) pp. 248--253.

\bibitem{Roberts:2007}C.A. Roberts, 
Recent results on blow-up and quenching for nonlinear Volterra equations. 
J. Comp. Appl. Math., 205 (2) (2007) pp. 736--743.

\bibitem{Smith:2011}H. Smith, 
An Introduction to Delay Differential Equations with Applications to the Life Sciences. 
Springer (2011)

\bibitem{Souplet:1998}P. Souplet, 
Blow-up in nonlocal reaction-diffusion equations. 
SIAM J. Math. Anal., 29 (6) (1998) pp. 1301--1334.

\bibitem{Straughan:2019}B. Straughan, 
Explosive Instabilities in Mechanics. Springer (2012) 

\bibitem{Takayasu:2017}A. Takayasu, K. Matsue, T. Sasaki, K. Tanaka, M. Mizuguchi, S. Oishi, 
Numerical validation of blow-up solutions of ordinary differential equations. 
J. Comp. Appl. Math., 314 (2017) pp. 10--29.

\bibitem{Yang:2013}Z. Yang, H. Brunner, 
Blow-up behavior of Hammerstein-type delay Volterra integral equations. 
Front. Math. China., 8 (2013) pp. 261--280.

\bibitem{Yang:2017}Z. Yang, T. Tang, J. Zhang, 
Blowup of Volterra integro-differential equations and applications to semi-linear Volterra diffusion equations. 
Num. Math., 10 (4) (2017) pp. 737--759. 

\bibitem{Zhou:2016}Y.C. Zhou, Z.W. Yang, H.Y. Zhang, Y. Wang, 
Theoretical analysis for blow-up behaviors of differential equations with piecewise constant arguments, 
Appl. Math. Comp., 274 (2016) pp, 353--361.
\end{thebibliography}

\end{document}